\documentclass[final]{siamltex}
\pdfoutput=1

\usepackage{graphicx}
\usepackage{amssymb,amsmath}
\usepackage{comment}
\usepackage{bbm,color}
\def\C{\mathbbm{C}}
\def\R{\mathbbm{R}}
\def\deg{d}         % degree of the polynomial filter
\def\kdim{m}        % Krylov subspace dimension (was previously \ell)
\def\numev{k}       % number of restart shifts per cycle
\def\rdim{r}        % number of restart roots
\def\CK{{\cal K}}   % calligraphic K for Krylov subspaces 
\def\CP{{\cal P}}   % calligraphic P for polynomials
\def\mrpoly{\pi}    % min res polynomial 
\def\CU{{\cal U}}   % calligraphic U
\def\CV{{\cal V}}   % calligraphic V
\def\Ugood{\CU}   % good invariant subspace
\def\Pgood{P_g}     % spectral projector for good invariant subspace
\def\Pbad{P_b}      % spectral projector for bad invariant subspace
\def\cost{\mbox{{\it cost}}}
\def\mvps{\mbox{{\it mvps}}}
\def\vops{\mbox{{\it vops}}}
\def\nev{\mbox{{\it nev}}}
\def\rtol{\mbox{{\it rtol}}}
\def\nnzr{\mbox{{\it nnzr}}}
\def\pof{\mbox{{\it pof}}@}
\def\MaxErr{\text{{\it MaxErr}}}
\def\MaxPof{\text{{\it MaxPof}}}

\mathcode`\@="8000
{\catcode`\@=\active \gdef@{\mkern1mu}} % See page 641 of Knuth's "Digital Typography"
\begin{document}

\title{Polynomial Preconditioned Arnoldi
\footnotemark[1] }

\author{Mark Embree\footnotemark[2]
\and Jennifer A. Loe\footnotemark[3]
\and Ronald B. Morgan\footnotemark[3] }

\maketitle

\renewcommand{\thefootnote}{\fnsymbol{footnote}}
\footnotetext[1]{The first author was supported by NSF grant DMS-1720257. The third author was supported by NSF grant DMS-1418677.}
\footnotetext[2]{Department of Mathematics and
Computational Modeling and Data Analytics Division,
Academy of Integrated Science,
Virginia Tech, Blacksburg, VA 24061
({\tt embree@vt.edu}).}
\footnotetext[3]{Department of Mathematics, Baylor
University, Waco, TX 76798-7328\\ ({\tt Jennifer\_Loe@baylor.edu}, {\tt Ronald\_Morgan@baylor.edu}).}
\renewcommand{\thefootnote}{\arabic{footnote}}

\begin{abstract}
Polynomial preconditioning can improve the convergence of the Arnoldi 
method for computing eigenvalues.  Such preconditioning significantly reduces the 
cost of orthogonalization; for difficult problems, it can also reduce the number of
matrix-vector products.  Parallel computations can particularly benefit from the
reduction of communication-intensive operations.
The GMRES algorithm provides a simple and effective way of generating the 
preconditioning polynomial.
For some problems high degree polynomials are especially effective, but they can lead to
stability problems that must be mitigated.  A two-level ``double polynomial preconditioning'' strategy provides an effective way to generate high-degree preconditioners.
\end{abstract}



\pagestyle{myheadings}
\thispagestyle{plain}
\markboth{M. EMBREE, J. A. LOE and R. B. MORGAN}{POLYNOMIAL PRECONDITIONED ARNOLDI}

\section{Introduction}

\begin{comment}
\begin{itemize}
\item \hemph{Do we want to emphasis more ``GMRES polynomial" or ``minres polynomial"?  In Abstract, we say GMRES.  Section 2 is titled with min res.  Algorithm in section 2 uses both. etc.  ME:  Let's go with ``GMRES polynomial.''
\end{itemize}
\end{comment}

We seek eigenvalues and eigenvectors of a large (possibly nonsymmetric) matrix $A$.  
The restarted Arnoldi algorithm~\cite{Saa80,So} (invoked by MATLAB's {\tt eigs} command) 
is a standard workhorse for such problems,  but for some matrices convergence is slow.  
On can improve convergence via a shift-invert transformation,
i.e., applying the algorithm to $(A-\mu I)^{-1}$ to find eigenvalues near $\mu\in\C$.\ \ 
Here we investigate an effective alternative that does not need any explicit inversion of $A$.
The \emph{polynomial preconditioned Arnoldi method}
 is fairly simple to implement and can accelerate convergence for difficult problems.

When applied to the matrix $A$ and starting vector $v$, the Arnoldi algorithm approximates eigenvalues using Rayleigh--Ritz estimates 
from the Krylov subspace 
\begin{equation} \label{eq:Kcond}
 \CK_\kdim(A,v) \equiv {\rm span}\{v, Av, \ldots, A^{\kdim-1}v\}.
\end{equation}
Any vector $x$ in this space, including the approximate eigenvectors,
can be written in the form $x = \omega(A) v$ for some 
$\omega\in \CP_{\kdim-1}$, where $\CP_s$ denotes the polynomials of
degree $s$ or less.
The Arnoldi process builds an orthonormal basis for the subspace~(\ref{eq:Kcond}) 
via a Gram--Schmidt process, requiring many inner products as $\kdim$ grows.

Polynomial preconditioning methods~\cite{La52B,LiXiVeYaSa,Rutis,Sa84b,Sa87b,Sa11,Sti58} apply the Arnoldi algorithm to the matrix $\pi(A)$, for some polynomial $\pi\in\CP_\deg$. Now eigenvalue estimates are drawn from the Krylov subspace
\begin{equation} \label{eq:preK}
 \CK_\kdim(\pi(A),v) = {\rm span}\{v, \pi(A)v, \ldots, \pi(A)^{\kdim-1}v\},
\end{equation}
a subspace of $\CK_{\deg(\kdim-1)+1}(A,v)$.
The large subspace  $\CK_{\deg(\kdim-1)+1}(A,v)$ 
contains better approximations to
the desired eigenvectors of $A$ than does $\CK_\kdim(A,v)$.
A polynomial preconditioner $\pi$ is effective if the
low-dimensional subspace $\CK_\kdim(\pi(A),v) \subseteq \CK_{\deg(\kdim-1)+1}(A,v)$ 
contains such improved estimates of the desired eigenvectors.

More specifically, 
any $x \in  \CK_\kdim(\pi(A),v)$ can be written as $x = \omega(\pi(A)) v$, 
where $\omega\in\CP_{\kdim-1}$.
Since $\omega\circ \pi \in \CP_{\deg(\kdim-1)}$, polynomial preconditioning leads to
eigenvector estimates that are high-degree polynomials in $A$.
This feature comes at a cost, since $\pi(A)$ must be applied to a vector
each time the subspace dimension $\kdim$ is increased.  
In typical high-performance computing environments these matrix-vector products can 
be evaluated more efficiently than inner products, which require significant communication and synchronization.
Thus the subspace~(\ref{eq:preK}) can be constructed much more efficiently (in terms of both work and storage)
than building out a standard Krylov subspace~(\ref{eq:Kcond}) of dimension $\deg(\kdim-1)+1$.
In summary, polynomial preconditioning gives an efficient way to involve 
high-degree polynomials while controlling the dimension of the subspace 
and limiting the cost of orthogonalization.

What is a good choice for the preconditioning polynomial $\pi$?
Section~\ref{sec:mrpoly} describes one choice for $\pi$ that is 
inspired by the GMRES algorithm.
By the spectral mapping theorem, every eigenvalue $\lambda$ of $A$ 
is mapped to the eigenvalue $\pi(\lambda)$ of $\pi(A)$.
As the convergence theory in section~\ref{sec:theory} illustrates,
effective preconditioners separate the desired eigenvalues
from the undesired ones.

Polynomial preconditioning is a special kind of \emph{spectral transformation}, in which the Arnoldi algorithm is applied to $f(A)$ for some function $f$ that maps the desired eigenvalues of $A$ to the largest magnitude eigenvalues of $f(A)$; see, e.g., \cite{MSR94}.  Typically such transformations involve a matrix inverse.  
One might seek a polynomial preconditioner $\pi$ that 
mimics a more complicated $f(A)$.  
For example, Thornquist~\cite{Tho06} advocates 
$\pi(A) \approx (A-\mu I)^{-1}$. 
Here we do not seek $\pi$ that approximates $(A-\mu I)^{-1}$, 
but merely one that distances the desired eigenvalues from the rest of the spectrum.

Although methods for polynomial preconditioning have been proposed in the past, they are not generally used in practice.  
\emph{To become popular, a polynomial preconditioner must be both effective 
and easy to implement.}
We shall also explore stability, an important consideration for practical algorithms.

Section~2 discusses the choice of the GMRES (minimum residual) polynomial for the preconditioning, followed by some convergence theory in Section~3. Numerical experiments begin in Section~4, and suggest several practical issues that a robust algorithm should address.  Section~5 studies the sensitivity of $\pi$ to the choice of GMRES starting vector, while Section~6 shows how to adjust the starting vector to make it more likely for Arnoldi to find the desired eigenvalues.  
Section~7 addresses numerical stability, suggesting the addition of duplicate roots in $\pi$ to better control distant unwanted eigenvalues. Finally, Section~8 describes \emph{double polynomial preconditioning}, which enables the use of very large degree polynomials.

%%%%%%%%%%%%%%%%%%%%%%%%%%%%%%%%%%%%%%%%%%%%%%%%%%%%%%%%%%%%%%%%%%%%%%%%%%%%%%%%
\section{Minimum Residual Polynomials} 
\label{sec:mrpoly}
%%%%%%%%%%%%%%%%%%%%%%%%%%%%%%%%%%%%%%%%%%%%%%%%%%%%%%%%%%%%%%%%%%%%%%%%%%%%%%%%

\begin{figure}[t!]
\includegraphics[scale=0.38]{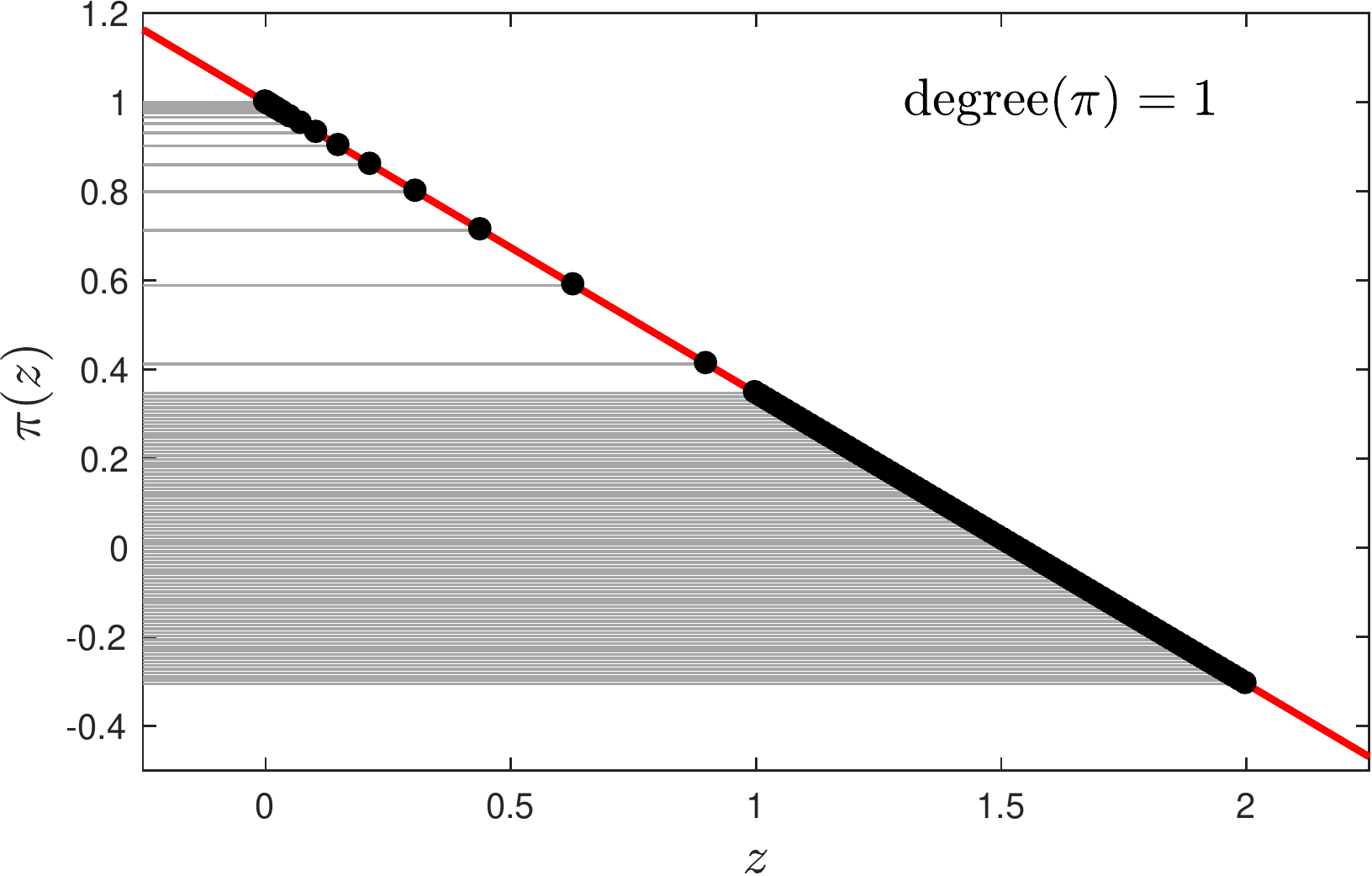}\quad
\includegraphics[scale=0.38]{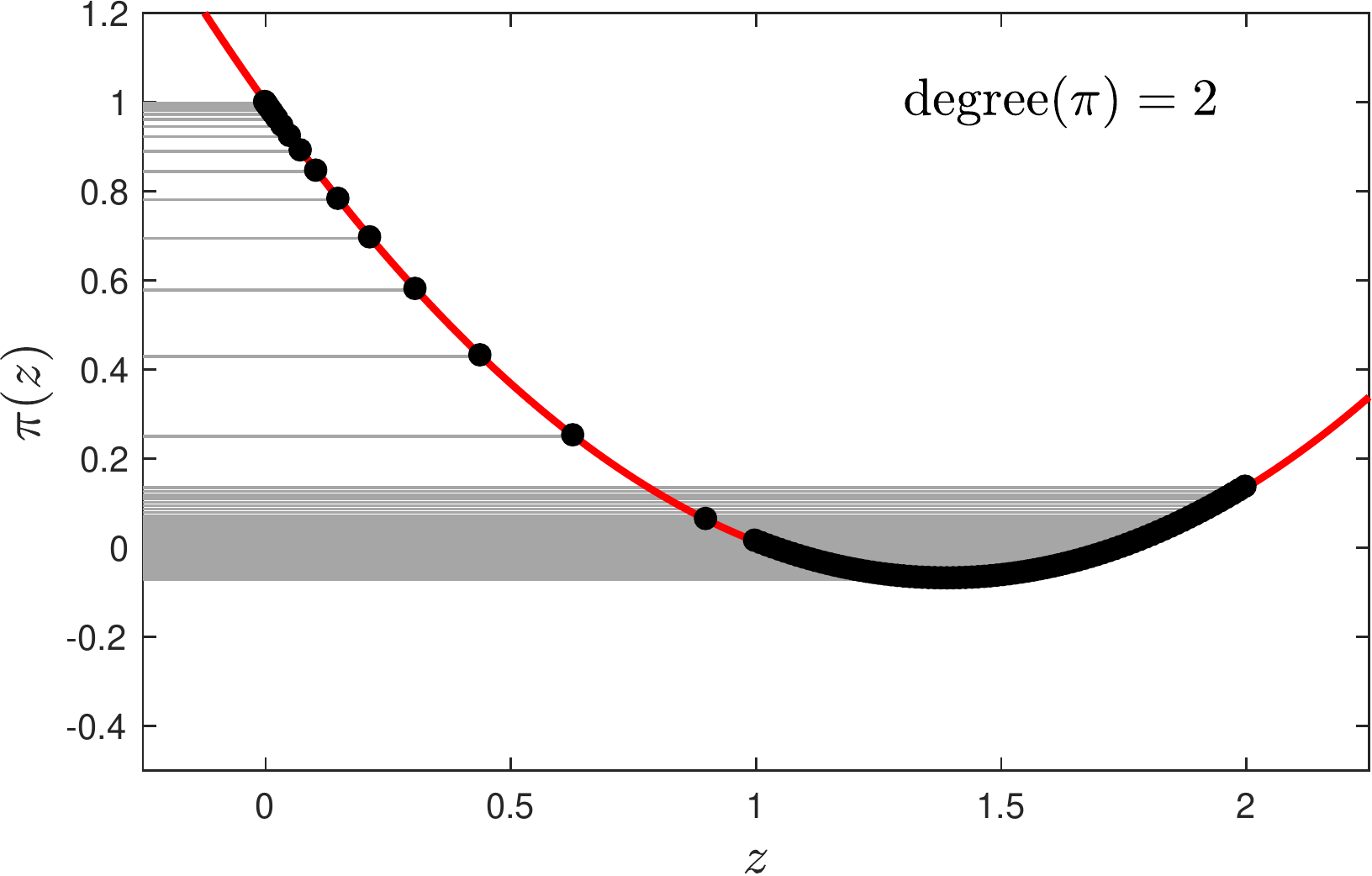}

\vspace*{.5em}
\includegraphics[scale=0.38]{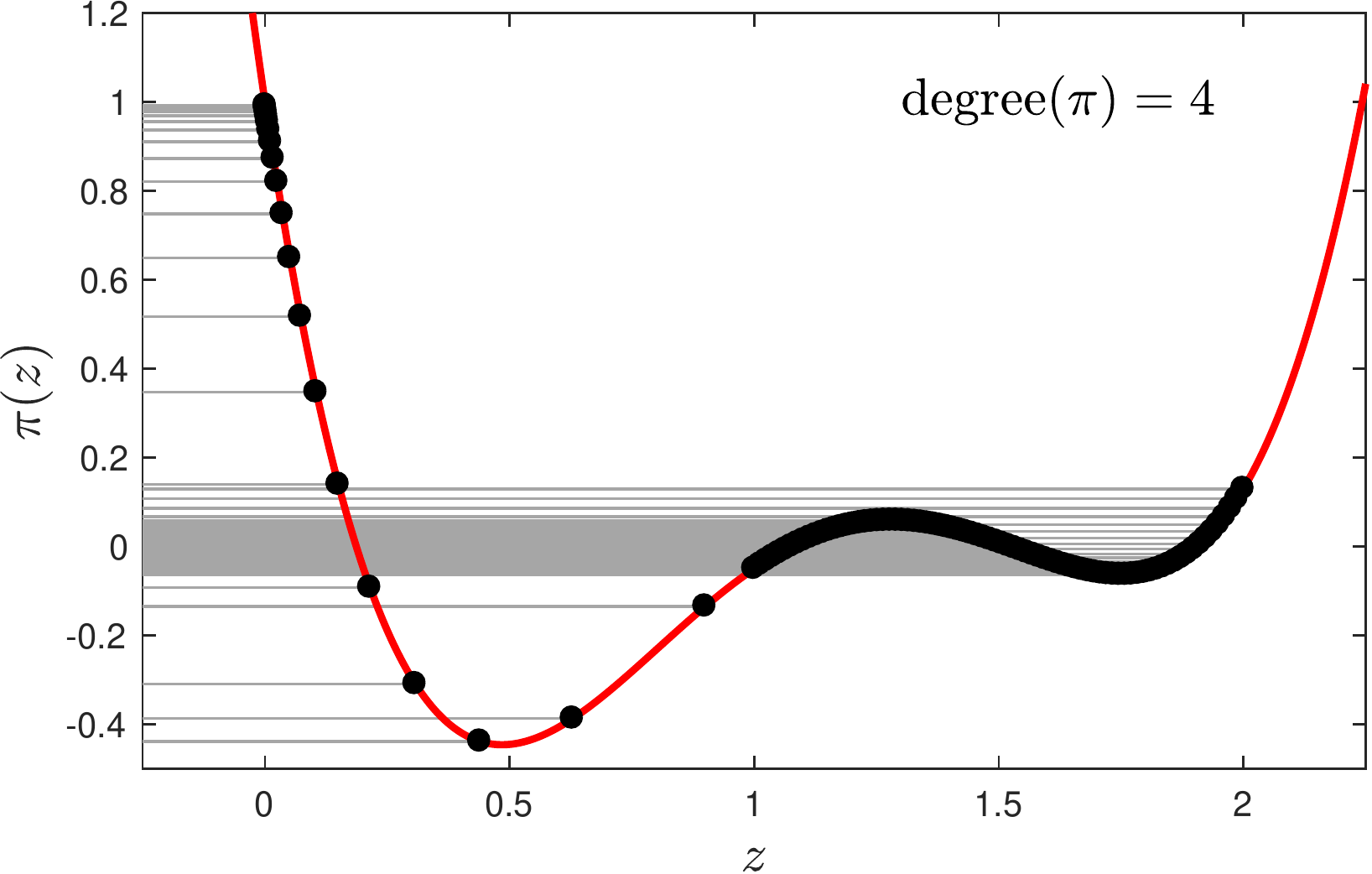}\quad
\includegraphics[scale=0.38]{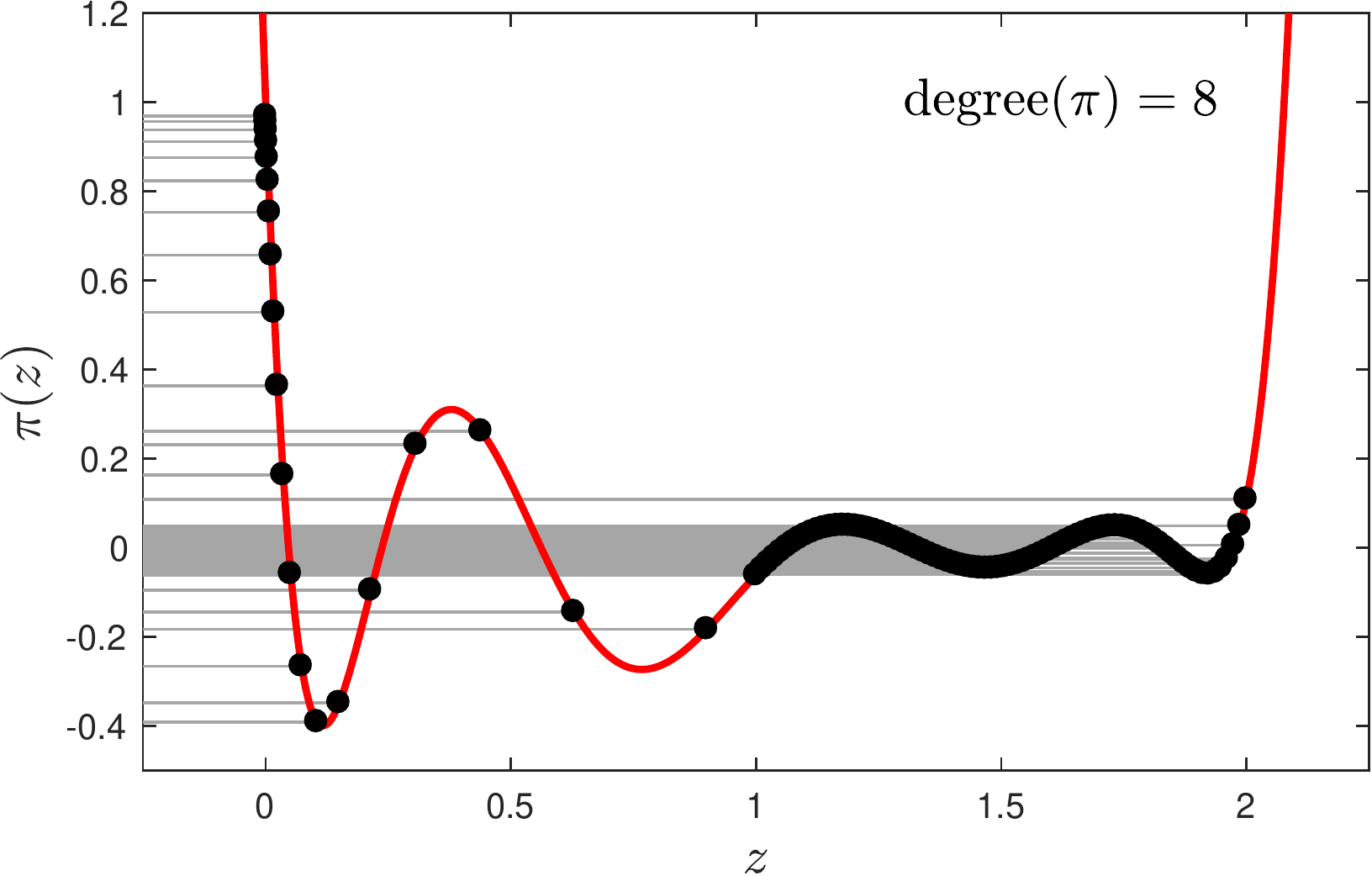}

\vspace*{-7pt}
\caption{\label{fig:polyplot}
GMRES polynomials (red lines) tend to separate eigenvalues
closest to the origin, while mapping large-magnitude eigenvalues of $A$ 
close to zero.  The black dots and horizontal gray lines show the 
values of $\mrpoly(\lambda_j)$.
}
\end{figure}

We seek the eigenvalues of $A$  nearest the origin.%
\footnote{If one seeks eigenvalues near $\mu \in \C$, 
replace $A$ with $A-\mu I$. If $\mu$ is in the interior of the spectrum, 
note that $|\pi(z)|$ is less likely to attain its maximum over the spectrum at $\mu$.  
The challenge of computing interior eigenvalues arises in Example~2, and is the subject of future work.}
For the polynomial preconditioner $\mrpoly$ we use the minimum residual polynomial~\cite{PPG,NaReTr} that arises when solving the linear system $Ax = b$ using the GMRES algorithm~\cite{SaSc} (or MINRES for symmetric $A$~\cite{PaSa}).  This polynomial satisfies
\[ \|\mrpoly(A) b\|_2 = \min_{\stackrel{\scriptstyle{p\in \CP_\deg}}
                            {\scriptstyle{p(0)=1}}}
                    \|p(A) b\|_2,\]
and hence $\mrpoly\in\CP_\deg$ must satisfy $\pi(0)=1$; $|\pi(z)|$ will generally 
be small over the spectrum of $A$.%
\footnote{This form allows one to write $\mrpoly(z) = 1-z \varphi(z)$ for some $\varphi\in \CP_{\deg-1}$.  Then $0 \approx \mrpoly(A) b = (I-A \varphi(A))b$ suggests that $\varphi(A) \approx A^{-1}$.  Thornquist uses this $\varphi$ as the polynomial preconditioner, replacing $A$ with $A-\mu B$ for generalized eigenvalue problems~\cite{Tho06}.}
% \hot{We refer to the GMRES polynomial as $q\in \CP_\deg$.}
Denote the eigenvalues of $A$ as $\sigma(A) = \{\lambda_j\}$, so the eigenvalues of $\mrpoly(A)$ are $\mrpoly(\lambda_j)$.
If GMRES converges quickly, then $\|\mrpoly(A)b\|$ is small, and the eigenvalues of $\mrpoly(A)$ are typically concentrated near~0.  However, the condition $\pi(0)=1$ means that $\mrpoly$ is generally not able to map the small eigenvalues of $A$ as close to zero,  making these eigenvalues better separated in the spectrum of $\mrpoly(A)$.  
Figure~\ref{fig:polyplot} illustrates this point,
using a symmetric $A$ with 20~eigenvalues logarithmically spaced in the interval $[10^{-3},0.9]$ 
and 80~eigenvalues uniformly spaced in the interval $[1,2]$; 
we seek a few of the smallest eigenvalues.  
Figure~\ref{fig:polyplot} shows $\mrpoly(z)$ for degree $\deg=1, 2, 4$ and $8$ 
(red lines) and the values of $\mrpoly(\lambda_j)$ (black dots and gray lines).
Since the small eigenvalues of $A$ are near the origin (where $\mrpoly(0)=1$),
$\mrpoly(\lambda_j)$ is large for these values, while $\mrpoly(\lambda_j)$
is small for the larger eigenvalues.  
Moreover, \emph{$\mrpoly$ separates the tightly clustered eigenvalues 
near the origin}.  The desired eigenvalues of $\mrpoly(A)$ (nearest~1) will be
easier for Arnoldi to compute than the corresponding (smallest) eigenvalues of $A$.
Figure~\ref{fig:polyplot} also hints at a complication: when the degree is large, the map $\mrpoly$ entangles some of the larger eigenvalues from the interval $[10^{-3},0.9]$ with those from $[1,2]$. 

%%%%%%%%%%%%%%%%%%%%%%%%%%%%%%%%%%%%%%%%%%%%%%%%%%%%%%%%%%%%%%%%%%%%%%%%%%%%%%%%
%%%%%%%%%%%%%%%%%%%%%%%%%%%%%%%%%%%%%%%%%%%%%%%%%%%%%%%%%%%%%%%%%%%%%%%%%%%%%%%%

The GMRES polynomial $\mrpoly$ is easy to construct in factored form.
To find its roots, run a cycle of GMRES($\deg$) and compute the harmonic Ritz values~\cite{Fr92,IE,PaPavdV} (reciprocals of Rayleigh--Ritz eigenvalue estimates for $A^{-1}$ from $A\CK_\deg(A,b)$).  For numerical stability, label the roots using the modified Leja ordering~\cite[alg.~3.1]{BaHuRe}, giving
\begin{equation} \label{eq:mrpoly}
\mrpoly(z) = \prod\limits_{i=1}^{\deg}\Big(1-\frac{z}{\theta_i}\Big).
\end{equation}
A quick listing of the algorithm follows.

%\vspace{.07in}
\newpage
\begin{center}
\textbf{Polynomial Preconditioned Arnoldi, \\ with GMRES polynomial of degree \boldmath $\deg$}
\end{center}
\medskip
\begin{enumerate}
\item {\bf Construction of the polynomial preconditioner, \boldmath$\mrpoly$:} 
\begin{enumerate}
\item Run one cycle of GMRES($\deg$).
\item Find the harmonic Ritz values, $\theta_1, \ldots, \theta_\deg$, which are the roots of the GMRES polynomial: with Arnoldi decomposition $AV_\deg = V_{\deg+1} H_{\deg+1,\deg}$, find the eigenvalues of $H_{\deg,\deg}^{} + h_{\deg+1,\deg}^2 @@f e_\deg^T$, where $f = H_\deg^{-*} e_\deg^{}$ with elementary coordinate vector $e_\deg=[0,\ldots,0,1]^T$.     
\item Order the GMRES roots with modified Leja ordering~\cite[alg.~3.1]{BaHuRe}.  (To avoid overflow or underflow for high degree polynomials, one can replace products of absolute values of differences of roots with addition of logarithms of these quantities.) 
\end{enumerate}
\medskip
\item {\bf PP-Arnoldi:} Apply restarted Arnoldi to the matrix $\pi(A) = \Pi_{i=1}^{d}(I - A/\theta_i)$.
(The experiments in Sections 4--8 use a thick-restarted Arnoldi method with exact shifts~\cite{Arnoldi-R,HRAM,StSaWu,WuSi}.)
\end{enumerate}
\vspace{.08in}

\begin{comment}
\begin{itemize}
\item \hemph{Incorporate example here, presently in Section~4?}
\end{itemize}
\end{comment}

%%%%%%%%%%%%%%%%%%%%%%%%%%%%%%%%%%%%%%%%%%%%%%%%%%%%%%%%%%%%%%%%%%%%%%%%%%%%%%%%
\section{Convergence theory for polynomial preconditioning}
\label{sec:theory}
%%%%%%%%%%%%%%%%%%%%%%%%%%%%%%%%%%%%%%%%%%%%%%%%%%%%%%%%%%%%%%%%%%%%%%%%%%%%%%%%

Let $\sigma(A)$ denote the spectrum of $A$.  
We seek some subset $\Sigma := \{\lambda_1, \ldots, \lambda_\numev\}$ of $\sigma(A)$, typically those of smallest magnitude.)
The eigenvalues in $\Sigma$ should be listed with the their full multiplicity,
but we presume that the eigenvalues in $\Sigma$ are \emph{nonderogatory}
(i.e., there exists only one linearly independent eigenvector for each
distinct eigenvalue in $\Sigma$); this assumption is standard for 
Krylov subspace convergence theory, required because of the concept 
of \emph{reachable invariant subspaces}~\cite{BER04,BES05}.

How does the preconditioned Arnoldi algorithm converge
as the dimension of the Krylov subspace increases?
Especially for nonsymmetric $A$, the error in the eigenvalue estimates can converge quite irregularly.  
Thus we prefer to analyze the rate at which the Krylov subspace in the
Arnoldi method ``captures'' the invariant subspace (span of eigenvectors 
and generalized eigenvectors) associated with $\Sigma$.

Let us make this notion precise.  
Let $\CU$ denote the maximal invariant subspace associated with the
desired eigenvalues $\Sigma$.\ \ Since none of these eigenvalues are 
derogatory, ${\rm dim}(\CU) = \numev$.
The Arnoldi algorithm approximates $\CU$ by some 
Krylov subspace $\CV$, whose dimension will generally differ from $\numev$.
The \emph{containment gap} (or just \emph{gap}) between $\CU$ and $\CV$,
\begin{equation} \label{def:gap}
  \delta(\CU,\CV) \equiv \max_{u \in \CU}\ \min_{v\in\CV}\ \frac{\|u-v\|}{\|u\|},
\end{equation}
measures the sine of the largest canonical angle between $\CU$ and its best
$\numev$-dimensional approximation from $\CV$.  
Note that $\delta(\CU,\CV) \in [0,1]$, with 
$\delta(\CU,\CV)=1$ when ${\rm dim}(\CV)<{\rm dim}(\CU) = \numev$:
$\CV$ must have
dimension at least $\numev$ in order to approximate all of $\CU$.
For additional properties of the gap, see~\cite{BER04,Cha83,Kat80}.  

For the polynomially preconditioned Arnoldi algorithm, 
the Krylov subspace $\CK_\kdim(\pi(A),v)$ plays the role of $\CV$.
How does $\delta(\CU,\CK_\kdim(\pi(A),v))$ depend on 
the choice of $\pi$, the dimension $\kdim$, and the starting vector $v$?
We shall apply the convergence theory from~\cite{BES05},
which uses the following notation and assumptions. 

\medskip
\begin{enumerate}
\item[(a)] Label the desired eigenvalues as $\lambda_1, \ldots, \lambda_\numev$,
       allowing multiplicities.  Assume none of these eigenvalues are 
       derogatory, and these eigenvalues are disjoint from the rest of
       the spectrum $\{\lambda_{\numev+1},\ldots, \lambda_n\}$.

 \item[(b)] $\Ugood$ denotes the $m$-dimensional invariant subspace associated 
      with $\lambda_1,\ldots, \lambda_\numev$. 

 \item[(c)] $\Pgood$ denotes the spectral projector onto $\Ugood$, so that
       $\Pbad \equiv I - \Pgood$ is the spectral projector onto the 
       invariant subspace associated with $\lambda_{\numev+1},\ldots, \lambda_n$.
 \item[(d)] $\Omega_b$ is a compact subset of $\C$ that contains the 
       eigenvalues $\lambda_{\numev+1}, \ldots, \lambda_n$.
 \item[(e)] $\alpha_g(z) \equiv (z-\lambda_1)\cdots (z-\lambda_\numev)$ 
       is the component of the minimal polynomial of $A$ associated with 
       the desired eigenvalues.
% \item[(f)] Assume $\Phi$ has $M$ distinct roots, and let $\Psi\in\CP_{M-1}$ be 
%       the polynomial that interpolates $1/\alpha_g$ at these roots.
% \item[(g)] Assume $\dim(\CU)  = \numev$ 
%             (the degree of the restart polynomial at each cycle), and that
%            $\kdim = 2@\numev$, a common setting for the implicitly 
%            restarted Arnoldi algorithm~\cite{LSY98,So}.
\end{enumerate}
For the sake of comparison, we start with Theorem~3.3 of~\cite{BES05}, which applies to the
standard Arnoldi method for $(A,v)$.
For a Krylov subspace of dimension
$\kdim \ge 2\numev$,
the theorem gives
\begin{equation} \label{eq:gapconv}
 \delta(\Ugood, \CK_\kdim(A, v)) \le 
   \left(\max_{\psi \in \CP_{\numev-1}} 
              \frac{\|\psi(A) \Pbad v\|}{\|\psi(A) \Pgood v\|}\right)
   \kappa(\Omega_b)
   \min_{p \in \CP_{\kdim-2@\numev}}  
   \max_{z\in \Omega_b} |1 - \alpha_g(z)@p(z)|.
\end{equation}
This bound has three ingredients.
\begin{itemize}
\item The starting vector $v$ only appears in the constant 
\[ C_1 \equiv \max_{\psi \in \CP_{\numev-1}} 
              {\|\psi(A) \Pbad v\| \over \|\psi(A) \Pgood v\|}.\]
If $\numev=1$, $C_1 = \|\Pbad v\|/\|\Pgood v\|$ gauges the bias of 
$v$ toward the desired invariant subspace.  For $\numev>1$, 
the eigenvalues of $A$ also influence $C_1$.
For additional details about $C_1$, see~\cite[sect.~5.1]{BER04}.

\item  The constant 
\[ C_2 \equiv \kappa(\Omega_b)\ge 1 \]
is a measure of the nonnormality of $A$ associated with the unwanted eigenvalues.
If $A$ is symmetric, then $\kappa(\Omega_b)=1$.
For nonnormal $A$, enlarging $\Omega_b$ to contain points beyond $\lambda_{\numev+1},\ldots, \lambda_n$ generally reduces $\kappa(\Omega_b)$.
For a detailed discussion of this constant, see~\cite[sect.~5.2]{BER04}.

\item As the Krylov subspace dimension $\kdim$ grows, 
the approximation problem 
\begin{equation} \label{eq:approx}
    \min_{p \in \CP_{\kdim-2@\numev}} 
   \max_{z\in \Omega_b} |1 - \alpha_g(z)@p(z)|.
\end{equation}
gives the mechanism for convergence.
The minimization problem seeks polynomials that approximate
$1/\alpha_g(z) = (z-\lambda_1)^{-1} \cdots (z-\lambda_\numev)^{-1}$ 
over $z\in\Omega_b$.  The convergence of this polynomial
approximation problem depends on the proximity of $\Omega_b$ 
to the desired eigenvalues $\lambda_1,\ldots, \lambda_\numev$.
Better separation between the desired and undesired eigenvalues 
yields faster convergence.
\end{itemize}

% \medskip
% Consider the special case of symmetric $A$ with 
% \[ \lambda_1 \le \cdots \le \lambda_m < \lambda_{m+1} \le \cdots \le \lambda_n.\]
% The set $\Omega_b = [\lambda_{m+1},\lambda_n]$; the symmetry of $A$ implies $C_2 = \kappa(\Omega_b) = 1$.  
% Note that $\lambda_m$ is the desired eigenvalue closest to $\Omega_b$, 
% and define
% \begin{equation} \label{eq:kappa}
%    \kappa := {\displaystyle \max_{z \in \Omega_b}\  |\lambda_m-z| 
%             \over 
%              \displaystyle \min_{z \in \Omega_b}\ |\lambda_m - z|} 
%            = {\lambda_n - \lambda_m \over \lambda_{m+1}-\lambda_m} \ge 1.
% \end{equation}
% Then using standard Chebyshev approximation theory~\hot{[XX]}, 
% the problem~(\ref{eq:approx}) converges at the
% asymptotic rate of 
% \begin{equation} \label{eq:convrate}
%     \rho := {\sqrt{\kappa} - 1 \over \sqrt{\kappa}+1} \in [0,1),
% \end{equation}
% meaning that~(\ref{eq:approx}) behaves like $C \rho^{\ell}$ 
% for some constant $C>0$ as $\ell\to\infty$.

%%%%%%%%%%%%%%%%%%%%%%%%%%%%%%%%%%%%%%%%%%%%%%%%%%%%%%%%%%%%%%%%%%%%%%%%%%%%%%%%
\subsection{Convergence theory for polynomial preconditioning}
\label{sec:polythy}
%%%%%%%%%%%%%%%%%%%%%%%%%%%%%%%%%%%%%%%%%%%%%%%%%%%%%%%%%%%%%%%%%%%%%%%%%%%%%%%%

Polynomial preconditioning alters the convergence bound~(\ref{eq:gapconv}),
replacing $A$ by $\pi(A)$.

%%%%%%%%%%%%%%%%%%%%%%%%%%%%%%%%%%%%%%%%%%%%%%%%%%%%%%%%%%%%%%%%%%%%%%%%%%%%%%%%
\begin{theorem} \label{thm:piconv}
Using the notation established above,
the gap between the desired invariant subspace $\Ugood$ 
(associated with the non-derogatory eigenvalues $\lambda_1,\ldots,\lambda_m$ of $A$)
and the polynomially-preconditioned Krylov subspace $\CK_\kdim(\pi(A),v)$ 
satisfies 
\[ \delta(\Ugood, \CK_\kdim(\pi(A), v)) 
 \le 
   \left(\max_{\psi \in \CP_{\numev-1}} 
              {\|\psi(\pi(A)) \Pbad v\| \over \|\psi(\pi(A)) \Pgood v\|}\right)
   \kappa(\Omega_b^\pi)
   \min_{p \in \CP_{\kdim-2@\numev}}  
   \max_{z\in \Omega_b^\pi} |1 - \alpha_g^\pi(z)@p(z)|,\]
provided
$\{\pi(\lambda_1),\ldots,\pi(\lambda_m)\} \cap 
 \{\pi(\lambda_{m+1}),\ldots, \pi(\lambda_n)\} = \emptyset$.
Here $\Omega_b^\pi$ is a compact subset of $\C$ that 
contains $\pi(\lambda_{\numev+1}), \ldots, \pi(\lambda_n)$,
and
\[ \alpha_g^\pi(z) := \big(z - \pi(\lambda_1)\big)\cdots \big(z - \pi(\lambda_\numev)\big).\]
\end{theorem}
For this bound to converge, $\pi$ must map the desired eigenvalues outside $\Omega_b^\pi$: 
\[ \{ \pi(\lambda_j)\}_{j=1}^\numev\ \cap\ \Omega_b^\pi\ =\ \emptyset.\]
Since the spectral projectors are invariant under the transformation $A \mapsto \pi(A)$ (provided the desired and undesired eigenvalues remain disjoint under $\mrpoly$), $\Pbad$ and $\Pgood$ are the same for $A$ and $\pi(A)$.

\medskip
Let us inspect this bound in the simplest case of symmetric $A$ with
$\numev=1$.
Since $\numev=1$, the constants in~(\ref{eq:gapconv}) and Theorem~\ref{thm:piconv}
that involve $v$ are just $C_1 = \|\Pbad v\| / \|\Pgood v\|$.
Since $A$ is symmetric, $C_2 = \kappa(\Omega_b) = \kappa(\Omega_b^\pi) = 1$.
The only difference in the convergence bounds in~(\ref{eq:gapconv}) and
Theorem~\ref{thm:piconv} 
then comes from the polynomial approximation problems.  
As suggested in Figure~\ref{fig:polyplot}, the map $\pi$ can effectively
separate the desired eigenvalues from the rest of the spectrum, 
making it possible that
\[ 
\min_{p \in \CP_{\kdim-2@\numev}}
   \max_{z\in \Omega_b^\pi} |1 - \alpha_g^\pi(z)@p(z)|
\ \ll\ 
\min_{p \in \CP_{\kdim-2@\numev}}
   \max_{z\in \Omega_b} |1 - \alpha_g(z)@p(z)|.
\] 
Keep $A$ symmetric but allow general $\numev\ge 1$.
Presuming that $\pi$ has real coefficients (so $\pi(\sigma(A)) \subset \R$),
denote
\[ \Omega_b^\pi \equiv \Big[ 
\min_{j=\numev+1,\ldots,n} \pi(\lambda_j),\ \max_{j=\numev+1,\ldots,n} \pi(\lambda_j)
\Big] \subset \R.\]
Let $\lambda_* \in \{\lambda_1,\ldots, \lambda_\numev\}$ denote the desired
eigenvalue that is mapped closest to $\Omega_b^\pi$:
\[ {\rm dist}(\pi(\lambda_*), \Omega_b^\pi) 
     = \min_{1\le j\le \numev} {\rm dist}(\pi(\lambda_j),\Omega_b^\pi)
     = \min_{1\le j\le \numev}\ \min_{z\in\Omega_b^\pi} |\pi(\lambda_j)-z|.\]
Supposing that $\pi(\lambda_*) \not\in \Omega_b^\pi$, define
\begin{equation} \label{eq:K}
        K \equiv {{ \max_{z \in \Omega_b^\pi}\  |\pi(\lambda_*)-z|} 
             \over 
              {\min_{z \in \Omega_b^\pi}\ |\pi(\lambda_*) - z|}}
            \ge 1.
\end{equation}
Then using standard Chebyshev approximation theory~\cite[sect.~6.11]{Saa03}, 
the polynomial approximation problem in Theorem~\ref{thm:piconv}
converges at the asymptotic rate
\begin{equation} \label{eq:convrate}
    \rho \equiv {\sqrt{K} - 1 \over \sqrt{K}+1} \in [0,1),
\end{equation}
meaning that there exists some constant $C>0$ such that for all $\kdim \ge 2@\numev$,
\[ \min_{p\in\CP_{\kdim-2\numev}} \max_{z\in\Omega_b^\pi} |1-\alpha_g^\pi(z)p(z)|
\ \le\  C \rho^{\kdim}.\]
Table~\ref{tbl:convrate} compares this convergence rate $\rho$ for the
standard Arnoldi method (first row, $\pi(z)=z$) to the convergence obtained for 
minimum residual polynomials $\pi$ of increasing degree $\deg$ 
for the symmetric $A$ used in Figure~\ref{fig:polyplot}.
We seek the $\numev=5$ smallest magnitude eigenvalues of $A$.
(When $\deg=1$, the shift-invariance property of the Krylov subspace 
implies that polynomial preconditioning and standard Arnoldi generate
the same approximation subspace: $\CK_\kdim(\pi(A), v) = \CK_\kdim(A,v)$.)
Note that the fastest convergence rate \emph{per polynomial degree}, $\rho^{1/d}$, is obtained for the standard case (or $d=1$).  This is expected:  the preconditioned space $\CK_\kdim(\pi(A), v)$ only contains a small part of the much larger space $\CK_{\deg(\kdim-1)+1}(A,v)$.  However, \emph{the  preconditioned method obtains this convergence rate with much lower dimensional subspaces}.

%%%%%%%%%%%%%%%%%%%%%%%%%%%%%%%%%%%%%%%%%%%%%%%%%%%%%%%%%%%%%%%%%%%%%%%%%%%%%%%%
\begin{table}[t!]
\caption{\label{tbl:convrate}
Asymptotic convergence rates for various polynomial degrees $\deg$, for the
example in Figure~\ref{fig:polyplot} with $\numev=5$ desired eigenvalues.
The first row is for standard Arnoldi, $\pi(A) = A$.
The set $\Omega_g^\pi$ is the smallest interval that contains 
$\{\pi(\lambda_j)\}_{j=1}^\numev$, 
while $\Omega_b^\pi$ is the smallest interval that contains 
the map of the unwanted eigenvalues,
$\{\pi(\lambda_j)\}_{j=\numev+1}^n$. 
When $\Omega_b^\pi$ and $\Omega_g^\pi$ overlap, we set $\rho=1$.
}
{\small \begin{center}\begin{tabular}{r c c c c}
   $\deg$    &  $\Omega_g^\pi$ &  $\Omega_b^\pi$  &  $\rho$  & $\rho^{1/\deg}$ \\[3pt]
\hline \\[-8pt]
\llap{standar}d&  $[ 0.0010,  0.0042]$ & $[\phantom{-}0.0060,  2.0000]$  & $0.9416$ & \\ 
   1  & $[ 0.9973,  0.9993]$ & $[-0.3055,  0.9961]$  & $0.9416$ & $0.9416$ \\ 
   2  & $[ 0.9936,  0.9985]$ & $[-0.0705,  0.9908]$  & $0.9030$ & $0.9503$ \\ 
   3  & $[ 0.9840,  0.9962]$ & $[-0.2018,  0.9772]$  & $0.8589$ & $0.9506$ \\ 
   4  & $[ 0.9690,  0.9925]$ & $[-0.4384,  0.9557]$  & $0.8232$ & $0.9525$ \\ 
   5  & $[ 0.9529,  0.9886]$ & $[-0.4155,  0.9329]$  & $0.7845$ & $0.9526$ \\ 
   6  & $[ 0.9292,  0.9829]$ & $[-0.4102,  0.8994]$  & $0.7407$ & $0.9512$ \\ 
   7  & $[ 0.9063,  0.9772]$ & $[-0.4036,  0.8673]$  & $0.7057$ & $0.9514$ \\ 
   8  & $[ 0.8753,  0.9695]$ & $[-0.3909,  0.8240]$  & $0.6650$ & $0.9503$ \\ 
  16  & $[ 0.5563,  0.8830]$ & $[-0.3501,  0.4002]$  & $0.4134$ & $0.9463$ \\ 
  24  & $[ 0.1046,  0.7217]$ & $[-0.3096,  0.2434]$  & $1.0000$ & $1.0000$ \\ 
\end{tabular}\end{center}}
\end{table}
%%%%%%%%%%%%%%%%%%%%%%%%%%%%%%%%%%%%%%%%%%%%%%%%%%%%%%%%%%%%%%%%%%%%%%%%%%%%%%%%

%%%%%%%%%%%%%%%%%%%%%%%%%%%%%%%%%%%%%%%%%%%%%%%%%%%%%%%%%%%%%%%%%%%%%%%%%%%%%%%%
\subsection{The effect of restarting on convergence}
\label{sec:restarts}
%%%%%%%%%%%%%%%%%%%%%%%%%%%%%%%%%%%%%%%%%%%%%%%%%%%%%%%%%%%%%%%%%%%%%%%%%%%%%%%%

To obtain accurate eigenvalue estimates while limiting the dimension $\kdim$ 
of the Krylov subspace~(\ref{eq:K}),
the Arnoldi algorithm is \emph{restarted}~\cite{Saa80,So}.  
We begin by describing the standard Arnoldi algorithm with no preconditioning.
After taking $\kdim$ steps of the Arnoldi process with the matrix $A$ and
starting vector $v_0 \equiv v$ (the first \emph{cycle}), 
the method runs a fresh set of $\kdim$ steps with the same $A$ 
but a new starting vector, $v_1$.
In general, cycle $c+1$ takes $\kdim$ Arnoldi steps with $A$ and starting vector $v_c$. 
These new starting vectors are formed using 
polynomial restart methods~\cite{Saa80,Sa84b,So}:
given some positive integer $\rdim \le \kdim - \numev$,
such methods construct $v_{c} = \phi_c(A) v_{c-1}$,
where $\phi_c \in \CP_\rdim$.
Aggregate these starting vectors to get $v_c = \Phi_c(A) v$,
for $\Phi_c \equiv \phi_1 \cdots \phi_{c}\in \CP_{c@ \rdim}$, so that
cycle $c+1$ of the restarted Arnoldi method uses the Krylov subspace
\begin{equation} \label{eq:Kfilter}
 \CK_\kdim(A,\Phi_c(A)v) = {\rm span}\{\Phi_c(A)v, A\Phi_c(A)v, \ldots, A^{\kdim-1}\Phi_c(A)v\},
\end{equation}
generally an $\kdim$-dimensional subspace of $\CK_{c@\rdim+\kdim}(A,v)$. 
For any $x \in \CK_\kdim(A,\Phi_c(A)v)$ 
there exists some $\omega\in\CP_{\kdim-1}$ such that 
\begin{equation}  \label{eq:omegaArs}
 x = \omega(A)\Phi_c(A)v,
\end{equation}
with $\omega \cdot \Phi_c \in \CP_{c@\rdim + \kdim-1}$.
Like polynomial preconditioning, restarting with polynomial filters
gives access to elements in a high-degree Krylov subspace, while keeping
the overall subspace dimension low.  

The convergence behavior will depend
on the polynomial filters~\cite{BER04,BES05}.
While these filters can be constructed from Chebyshev polynomials~\cite{Sa84b}, filters built from ``exact shifts'' 
(unwanted Ritz values)~\cite{So} are simpler and more widely used.
Such shifts always give convergence for symmetric $A$~\cite{So} 
(aside from pathological $v_0$), 
and usually work well for nonsymmetric $A$ 
(though they can fail, in theory~\cite{DM12,Emb09}).
Moreover, these filters can be implemented stably using the
implicitly restarted Arnoldi~\cite{So}, thick-restart Arnoldi~\cite{Arnoldi-R,HRAM,StSaWu,WuSi} and Krylov--Schur~\cite{St01} algorithms.
Since these exact-shift filters are built up during each cycle of the 
restarted Arnoldi procedure, they require numerous inner product evaluations.

Practical computations will combine polynomial preconditioning with
polynomial restarting, using the approximation space
\begin{eqnarray} 
 \CK_\kdim(\pi(A),\Phi_c(\pi(A))v) &=& \nonumber \\
&& \hspace*{-5em}  {\rm span}\{\Phi_c(\pi(A))v, \pi(A)\Phi_c(\pi(A))v, \ldots, \pi(A)^{\kdim-1}\Phi_c(\pi(A))v\}.
\label{eq:Kfilterpre}
\end{eqnarray}
Generally~(\ref{eq:Kfilterpre}) is an $\kdim$-dimensional subspace of
$\CK_{\deg(c@\rdim+\kdim-1)}(A,v)$.
Now for any $x\in \CK_\kdim(\pi(A),\Phi_c(\pi(A))v)$ 
there exists a polynomial $\omega \in \CP_{\kdim-1}$ such that
\begin{equation} \label{eq:omegapiA}
   x = \omega(\pi(A))\Phi_c(\pi(A)) v,
\end{equation}
so preconditioning gives access to vectors $x$ that can be written 
in terms of a polynomial ($(\omega \circ \pi)\cdot(\Phi_c \circ \pi)$)
of degree $\deg@(c@\rdim+\kdim-1)$, i.e., 
$d$ times larger than available with restarting alone in~(\ref{eq:omegaArs}).

\begin{figure}[t]
\includegraphics[width=5in]{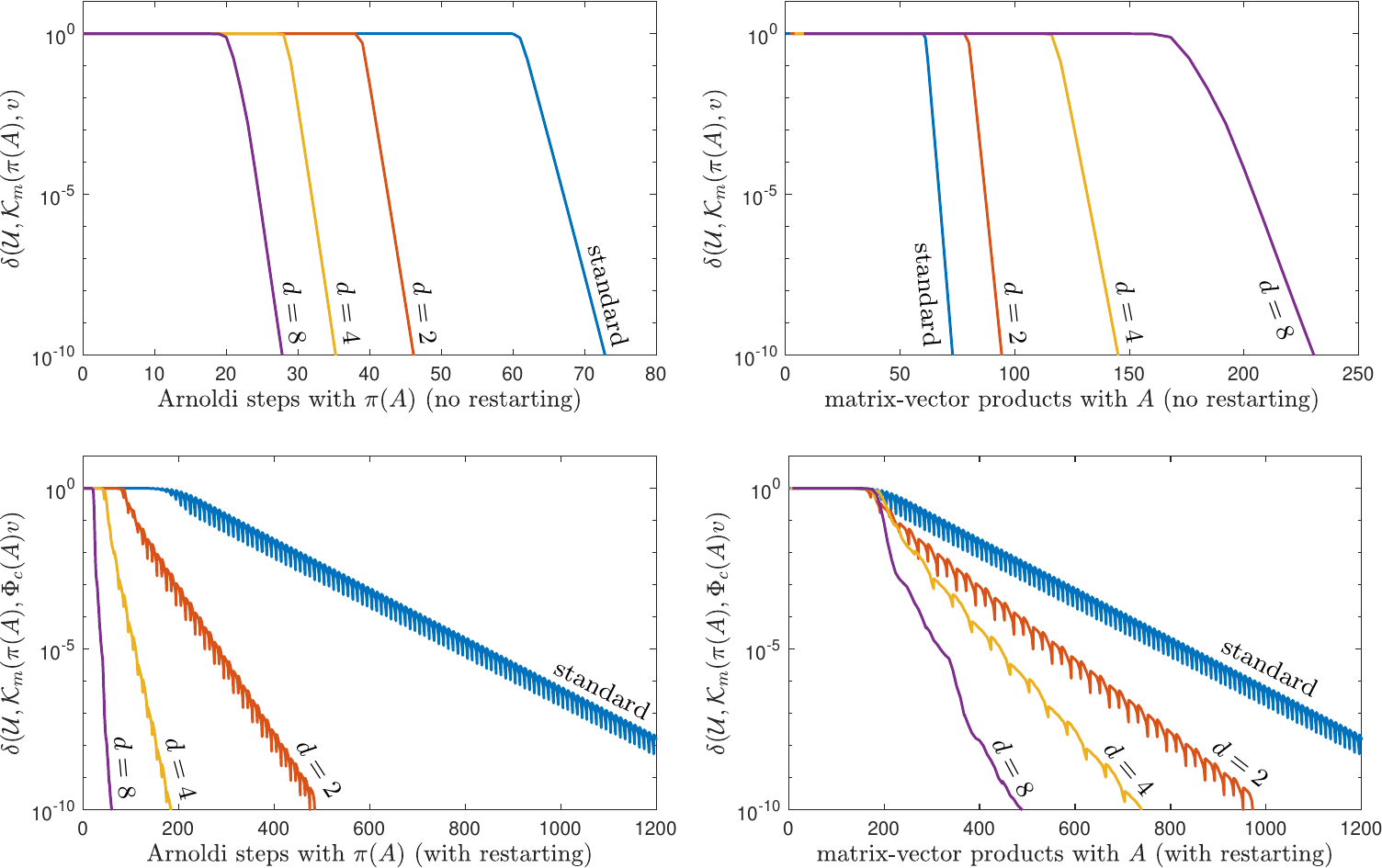}

\caption{\label{fig:prearn}
Arnoldi iterations and matrix-vector products for various choices of $d$
and $\numev=5$ eigenvalues,
for the example in Figure~\ref{fig:polyplot} without restarts (top) and
with restarts using $\kdim = 2@\numev$ (bottom).
Without restarts,  standard Arnoldi requires more steps
(top left) but fewer matrix-vector products (top right).  
Restarting can give polynomial preconditioning another advantage,
making it faster in terms of both steps
(bottom left) and matrix-vector products (bottom right).
}
\end{figure}

How do preconditioning and restarting combine to affect convergence?
We first address this question by continuing the experiment started
in Figure~\ref{fig:polyplot}, seeking the $\numev=5$ smallest magnitude eigenvalues.
Figure~\ref{fig:prearn} compares convergence of the polynomially-preconditioned
Arnoldi algorithm, with and without restarting.
In the top figures (no restarts), increasing the preconditioning
polynomial degree $\deg$ improves the convergence, but increases the
number of matrix-vector products involving $A$.
The bottom figures incorporate polynomial restarting, 
limiting the Krylov subspace to have dimension $\kdim=2@\numev$ 
and using filter polynomials of degree $\rdim = \numev$ at each cycle.
Polynomial preconditioning improves convergence by a larger margin,
enough to also deliver convergence using fewer matrix-vector products.

%%%%%%%%%%%%%%%%%%%%%%%%%%%%%%%%%%%%%%%%%%%%%%%%%%%%%%%%%%%%%%%%%%%%%%%%%%%%%%%%
%%%%%%%%%%%%%%%%%%%%%%%%%%%%%%%%%%%%%%%%%%%%%%%%%%%%%%%%%%%%%%%%%%%%%%%%%%%%%%%%

By adapting~\cite[eq.~(3.10)]{BES05}, 
we can also provide a bound on the gap between the 
desired invariant subspace $\Ugood$ and the restarted Krylov subspace
using polynomial preconditioning.
Suppose the aggregate polynomial filter $\Psi_c \in \CP_{c@\rdim}$
has $R$ distinct roots (the shifts), also distinct from the images of the desired eigenvalues,
$\Sigma^\mrpoly = \{\mrpoly(\lambda_1),\ldots, \mrpoly(\lambda_\numev)\}$.
Let $\Psi \in \CP_{R-1}$ interpolate $1/\alpha_g^\pi$ at these $R$ points.
Then
\begin{equation} \label{eq:gapconvrs}
\delta(\CU, \CK_\kdim(A, \Phi_c(A) v)) \le 
   \left(\kern-1pt \max_{\psi \in \CP_{\numev-1}} \kern-1pt
              {\|\psi(\pi(A)) P_b v\| \over \|\psi(\pi(A)) P_g v\|}\kern-1pt\right)\kern-1pt
   \kappa(\Omega_b^\pi) 
   \max_{z\in \Omega_b^\pi} |1 - \alpha_g^\pi(z)\Psi_c(z)|.
\end{equation}
If the roots of $\Psi_c$ are distributed throughout $\Omega_b^\pi$,
we expect $\max_{z\in\Omega_b^\pi} |1 - \alpha_g^\pi(z)\Psi_c(z)|$ to be small.
While descriptive bounds on the location of the shifts for 
 nonsymmetric $A$ (or $\mrpoly(A)$) have proved elusive~\cite{DM12,Emb09},
the bound~(\ref{eq:gapconvrs}) provides an indication of how 
the polynomial preconditioner and the shifts can combine to influence convergence.

%%%%%%%%%%%%%%%%%%%%%%%%%%%%%%%%%%%%%%%%%%%%%%%%%%%%%%%%%%%%%%%%%%%%%%%%%%%%%%%%
\section{Experiments} \label{sec:ex}
% * <ronald_morgan@baylor.edu> 2018-03-14T19:03:56.718Z:
% 
% Change name of section.
% 

How does polynomial preconditioning perform in practice?
In this section we explore three examples involving nonsymmetric $A$.
Here and in future sections, our experiments use the thick restarted Arnoldi method~\cite{HRAM,WuSi}, which we refer to as Arnoldi$(m,k)$: 
$m$ denotes the largest subspace dimension and $k<m$ denotes the number of Ritz values kept at each restart.%
\footnote{All examples involve real matrices; to preserve real arithmetic during the restarting process, sometimes $k$ is temporarily reduced to $k-1$ to avoid splitting a conjugate pair of Ritz values.}  
We seek the $\nev < k$ smallest magnitude eigenvalues, leaving a buffer of $k-\nev$ eigenvalues to accelerate convergence.
Each orthogonalization step is followed by a pass of reorthogonalization.  
Convergence is tested at each Arnoldi cycle using the original matrix $A$:
Let $\nu_1, \ldots, \nu_m$ denote the Ritz values for $\pi(A)$ at the end of a cycle, ordered by increasing distance from~1 ($|1-\nu_1| \le |1-\nu_2| \le \cdots \le |1-\nu_m|$), and let $y_1, \ldots, y_m\in\C^n$ denote the associated unit-norm Ritz vectors.  Then the Rayleigh quotient $\mu_j^{}\equiv y_j^*Ay_j^{}$ gives an eigenvalue estimate for $A$.
When seeking $\nev < k$ eigenvalues, we require $\|A y_j - \theta_j y_j\| \le \rtol$ for
$j=1,\ldots, \nev$, with $\rtol = \|A\|@10^{-8}$ unless otherwise stated.  
The reported operation counts give a rough impression of the matrix-vector products with~$A$,
dot products, and other vector operations (scalar multiplications and additions); these counts
will vary a bit with implementation details (e.g., reorthogonalization strategy, which residuals are checked at each cycle, etc.).  In each case, to explore robustness we use a random starting vector to generate $\mrpoly$, and a different random starting vector for the Arnoldi iterations.
(In practice one might naturally use the same starting vector to generate $\mrpoly$ and for the Arnoldi iterations.)  We average our results over ten trials, to minimize variation due to these starting vectors.

{\it Example 1.}
Consider a second-order finite difference discretization of a convection-diffusion equation on the unit square with homogeneous Dirichlet boundary conditions.  On the bottom half of the square, the operator is $ - u_{xx} - u_{yy} + 20u_{x} $; on the top half, $ - 100 u_{xx} - 100 u_{yy} + 2000 u_{x} $.  We use five increasingly fine discretizations that give matrices of size $n = 2500$, $\mbox{10,000}$, $\mbox{40,000}$, $\mbox{160,000}$ and $\mbox{640,000}$.
We seek the~$\nev=15$ smallest magnitude eigenvalues and the corresponding eigenvectors using Arnoldi(50,20), meaning the maximum subspace dimension is $m = 50$, and $k = 20$ Ritz vectors are saved at the restart.

%%%%%%%%%%%%%%%%%%%%%%%%%%%%%%%%%%%%%%%%%%%%%%%%%%%%%%%%%%%%%%%%%%%%%%%%%%%%%%%%
\begin{figure}
%\vspace{-2.6in}
%\hspace{-.7in}
%\includegraphics[scale=.75]{ppaf2a.pdf}
%\vspace{-2.8in}
\begin{center}
\includegraphics[scale=.63]{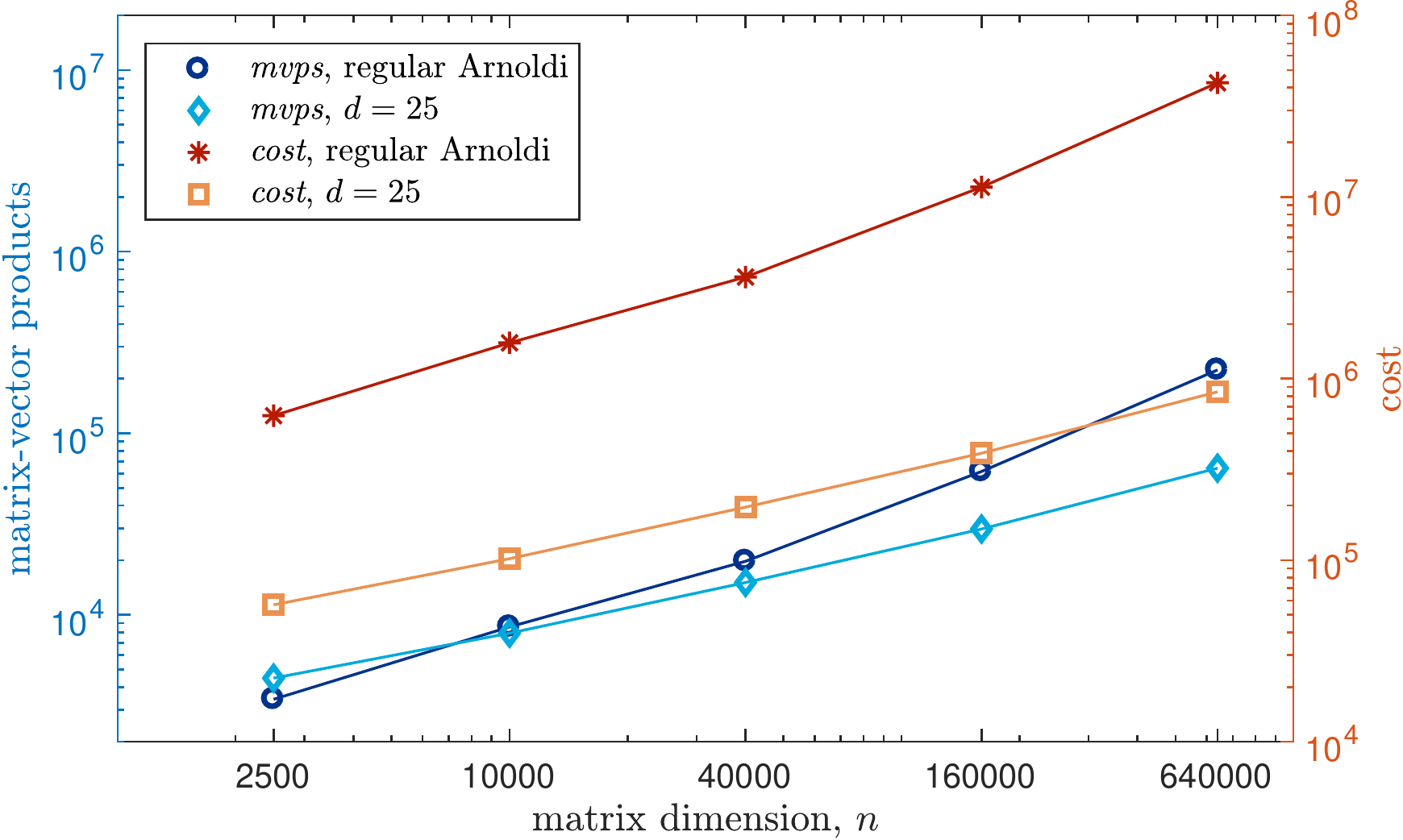}
\end{center}

\vspace*{-5pt}
\caption{\label{fig:ppaf2a}
Example 1 (convection diffusion matrix): Comparison of
regular Arnoldi(50,20) to polynomial preconditioned Arnoldi(50,20)
with degree $\deg=25$ for matrices of size $\mbox{2,500}$, $\mbox{10,000}$, $\mbox{40,000}$, $\mbox{160,000}$ and $\mbox{640,000}$.   Circles (regular Arnoldi) and diamonds (polynomial preconditioned) indicate the number of matrix-vector products (left axis).   Asterisks and squares show the corresponding approximate cost ($\cost = \nnzr\times \mvps + \vops$, right axis).}
\end{figure}
%%%%%%%%%%%%%%%%%%%%%%%%%%%%%%%%%%%%%%%%%%%%%%%%%%%%%%%%%%%%%%%%%%%%%%%%%%%%%%%%

Figure~\ref{fig:ppaf2a} compares regular Arnoldi(50,20) to degree $\deg=25$ polynomial preconditioned Arnoldi(50,20) for the five matrices, giving both the number of matrix-vector products for convergence (left axis) and an estimate of the total cost (right axis).  
% Matrix-vector products are associated with the left axis (circles for regular and diamonds for polynomial preconditioned).  
For $n = \mbox{2,500}$, polynomial preconditioning uses slightly more matrix-vector products.   As the matrix size increases and the eigenvalue problem becomes more difficult, polynomial preconditioning becomes increasingly better in comparison. 
For $n = \mbox{640,000}$, regular Arnoldi averages $\mbox{223,103}$ matrix-vector products, while polynomial preconditioned Arnoldi with $\deg=25$ only averages  $\mbox{64,227.3}$.  
% (Table~\ref{tbl:cdo} gives additional details about the operation counts for this problem.)

The computational cost is estimated as $\cost = \nnzr\times\mvps + \vops$, where 
$\nnzr \approx 5$ is the average number of nonzeros per row in $A$, 
$\mvps$ is the number of matrix-vector products, and $\vops$ is the number of length-$n$ vector operations, such as dot products and daxpy's.  Of course, the cost of an 
${\it mvp}$ compared to a ${\it vop}$ depends on the computer and implementation, but this estimate shows the potential for polynomial preconditioning to reduce computational cost.  
In Figure~\ref{fig:ppaf2a}, $\cost$ is associated with the right axis (asterisks for regular Arnoldi, squares for $\deg=25$).  The axes are scaled so the values on the left axis are one-fifth of the corresponding height on the right axis, allowing one to see what portion of $\cost$ is due to $\mvps$.  For regular Arnoldi, most matrix-vector products are accompanied by an orthogonalization step; preconditioning uses more matrix-vector products (to compute $\mrpoly(A)v$) relative to vector operations.  For this sparse $A$,  $\vops$  dominate  $\mvps$ for regular Arnoldi.  
Polynomial preconditioning with $\deg=25$ reduces the $\vops$ per matrix-vector product by a factor of about~22 for all the $n$ values shown here.  
With $n=\mbox{2,500}$, the $\cost$ for regular Arnoldi is $\mbox{624,768}$, but only $\mbox{56,755}$ for polynomial preconditioning.  With $n = \mbox{640,000}$, the comparison is $\mbox{42,534,104}$ to $\mbox{845,416}$.  Increasing the degree to $\deg = 100$ drops the cost to $\mbox{524,137}$: over 80~times cheaper than regular Arnoldi.

%        cost for regular, n=2500:     624767.95 
%           cost for d=25, n=2500:      56754.86 
%      cost for regular, n=640000:   42534104.48 
%         cost for d=25, n=640000:     845416.66 
%         cost for d=50, n=640000:     612165.61 
%        cost for d=100, n=640000:     524136.67 

\begin{figure}[t!]
%\vspace{-2.7in}
%\hspace{-.7in}
%\includegraphics[scale=0.75]{ppaf2pnew}
%\vspace{-2.8in}
\begin{center}
\includegraphics[width=3.75in]{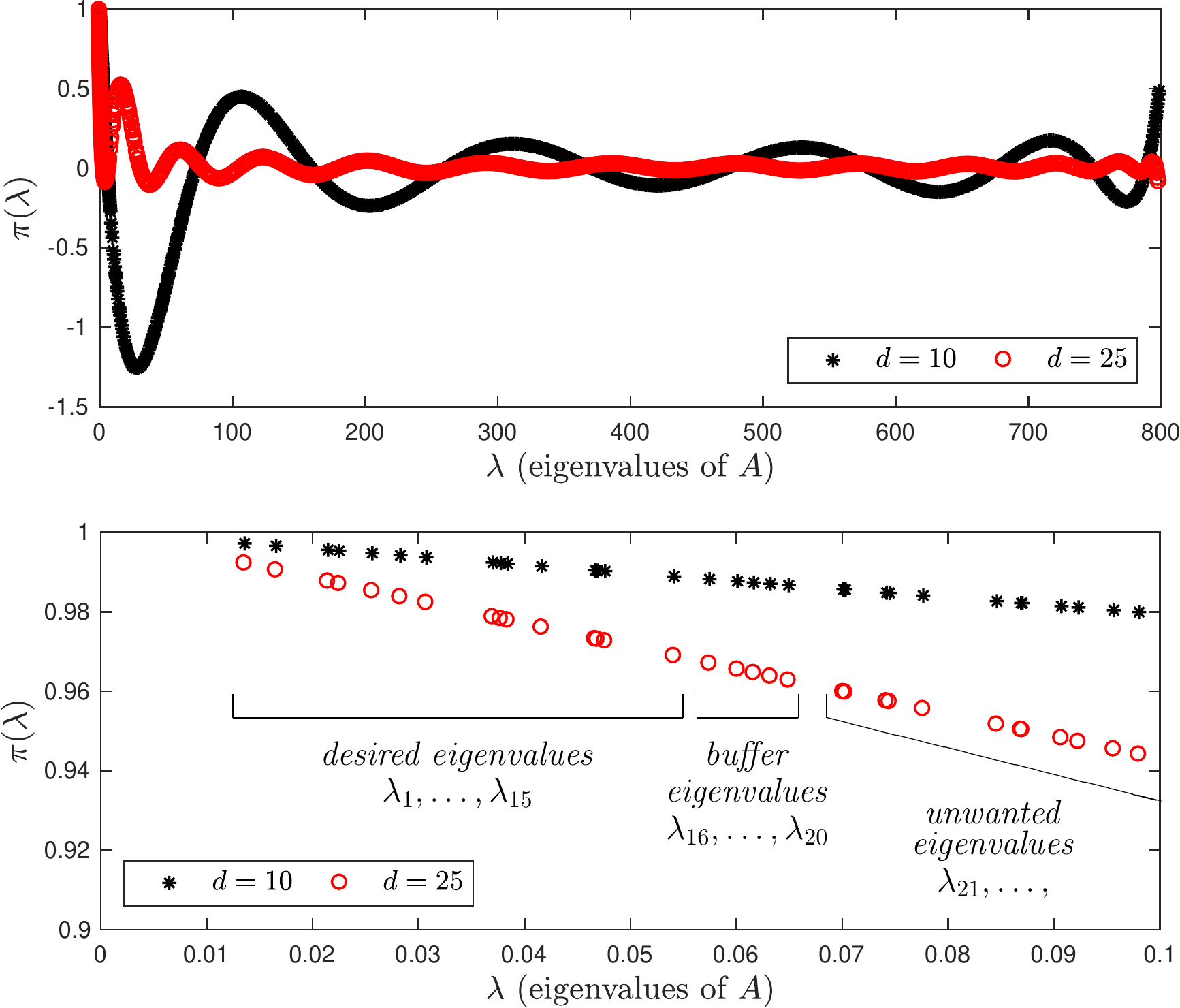}
\end{center}

\vspace*{-5pt}
\caption{\label{fig:ppaf2pnew}
Example 1 (convection diffusion matrix, $n=\mbox{10,000}$):
Polynomials of degree $\deg=10$ (black) and $\deg=25$ (red) for the convection-diffusion matrix of size $n=\mbox{10,000}$.  The upper plot shows the polynomial values $\pi(\lambda)$ plotted for all eigenvalues $\lambda$ of $A$.  The bottom plot zooms in on the smallest eigenvalues, highlighting the sought-after $\nev=15$ smallest magnitude eigenvalues, the buffer of $k-\nev=5$ additional eigenvalues, and a few of the remaining unwanted eigenvalues.} 
\end{figure}

We continue the example by looking at the polynomials for the matrix of size $\mbox{10,000}$.  The upper portion of Figure~\ref{fig:ppaf2pnew} shows representative GMRES polynomials of degree~10 and~25, evaluated at all of the eigenvalues.  The steeper slope at the origin of the $\deg=25$ polynomial better separates the small eigenvalues from the others, as is clear from the close-up plot on the bottom.  The desired first~$\nev=15$ eigenvalues come first, followed by the next $k-\nev = 5$ eigenvalues that serve as a buffer for the desired ones, and finally the next few eigenvalues:
the polynomial mapping separates the small eigenvalues from the others.  
A gap ratio for the 15th eigenvalue that takes into account the buffer Ritz values is 
$\frac{\lambda_{21}-\lambda_{15}}{\lambda_{10000}-\lambda_{21}} = 2.00\times10^{-5}$ 
for the matrix $A$, which gives some indication of the eventual convergence of that eigenvalue.  With $\deg = 10$ polynomial preconditioning, the gap ratio improves to
$\frac{\mrpoly(\lambda_{21})-\mrpoly(\lambda_{15})}{\mrpoly(\lambda_{5099})-\mrpoly(\lambda_{21})} = 1.14\times 10^{-3}$, 
and for $\deg=25$, it is still better: $\frac{\mrpoly(\lambda_{21})-\mrpoly(\lambda_{15})}{\mrpoly(\lambda_{5145})-\mrpoly(\lambda_{21})} = 1.62\times10^{-2}$.
While these ratios suggest that polynomial preconditioning improves the convergence rate in terms of Arnoldi cycles, a larger gap ratio \emph{does not} guarantee a reduction of overall matrix-vector products, since each iteration with a higher degree polynomial requires more of them.  As Figure~\ref{fig:ppaf2a} shows, for $n=\mbox{10,000}$
the number of matrix-vector products for $\deg=25$ is a bit smaller than for regular Arnoldi.
%For this example, matrix-vector products with the degree 25 are reduced (by 29\%), as was seen in Figure~\ref{fig:ppaf2a}.

%  gap ratio (d=0): 2.0032988e-05
% gap ratio (d=10): 1.1369308e-03 (min value at index = 5099)
%  gap ratio (d=0): 1.5627467e-05
% gap ratio (d=25): 1.6243782e-02 (min value at index = 5145)

%%%%%%%%%%%%%%%%%%%%%%%%%%%%%%%%%%%%%%%%%%%%%%%%%%%%%%%%%%%%%%%%%%%%%%%%%%%%%%%%
%%%%%%%%%%%%%%%%%%%%%%%%%%%%%%%%%%%%%%%%%%%%%%%%%%%%%%%%%%%%%%%%%%%%%%%%%%%%%%%%

%%%%%%%%%%%%%%%%%%%%%%%%%%%%%%%%%%%%%%%%%%%%%%%%%%%%%%%%%%%%%%%%%%%%%%%%%%%%%%%%
\begin{figure}[h!]
% \vspace{-2.65in}
% \hspace{-.7in}
% \includegraphics[scale=.75]{ppaf2b}
% \vspace{-2.8in}
\begin{center}
\includegraphics[scale=.63]{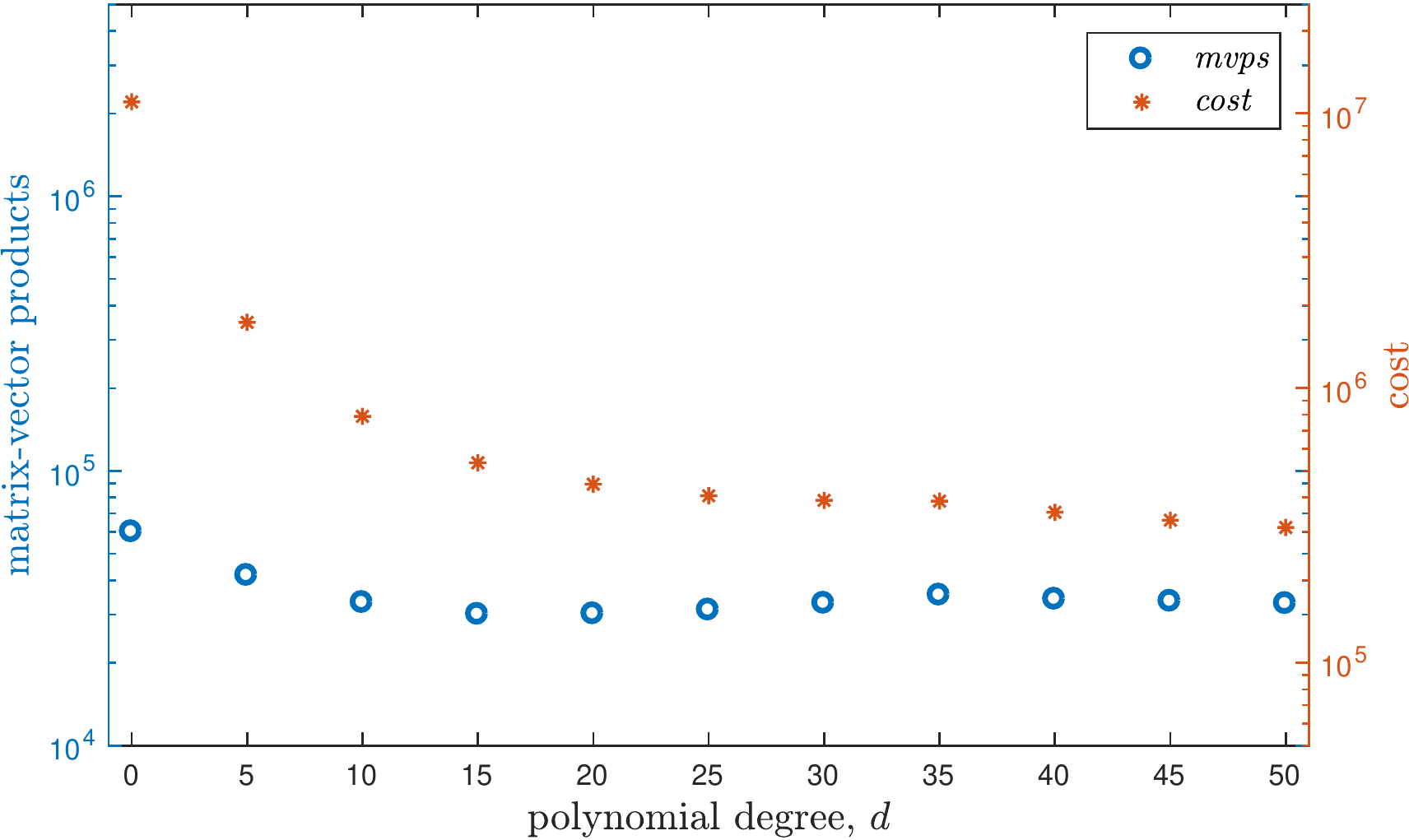}
\end{center}

\vspace*{-5pt}
\caption{\label{fig:ppaf2b}
Example 1 (convection diffusion matrix, $n=\mbox{160,000}$):  A comparison of convergence between standard Arnoldi(50,20) ($\deg=0$) and polynomial preconditioned Arnoldi(50,20) with $\deg = 5, 10, 15, \ldots, 50$.   Circles indicate the number of matrix-vector products (left axis).   The approximate cost ($\cost = \nnzr \times \mvps + \vops$) is indicated with asterisks (right axis).  
% \hemph{Figure~\ref{fig:ppaf2b} could be switched over to the larger matrix of size 640,000.  That is the size used in Table 8.1.}
}
\end{figure}
%%%%%%%%%%%%%%%%%%%%%%%%%%%%%%%%%%%%%%%%%%%%%%%%%%%%%%%%%%%%%%%%%%%%%%%%%%%%%%%%

Figure~\ref{fig:ppaf2b} shows results for the matrix of size $n = \mbox{160,000}$ using different degree polynomials.  Standard Arnoldi(50,20) is plotted at $\deg=0$; then polynomial preconditioned Arnoldi(50,20) is applied with $\deg = 5, 10, 15, \ldots, 50$.  The matrix-vector products, denoted on the left axis, hit a minimum for $\deg = 15$ and increase slightly beyond that.  The cost estimate (right axis) decreases.
For $\deg = 50$, $\cost \approx \mbox{316,060}$.  
The $\cost$ can go down a bit further with higher $\deg$; 
for $\deg = 150$, $\cost\approx  \mbox{298,308}$.
See Table~\ref{tbl:cdo} for further testing with a larger version of this example, including timings and dot product counts.

% minimum matvecs =   29976.50 attained at d = 15
% d =   0; cost = 11012783.79
% d =   5; cost = 1733894.62
% d =  10; cost =  788956.80
% d =  15; cost =  534862.44
% d =  20; cost =  447569.70
% d =  25; cost =  406649.54
% d =  30; cost =  391981.02
% d =  35; cost =  390335.30
% d =  40; cost =  356294.37
% d =  45; cost =  334718.05
% d =  50; cost =  316059.69
% d = 150; cost =  298307.64
%
% For example 1, the desired eigenvalues were always found to mild accuracy: there was no problem converging to incorrect eigenvalues.

\begin{figure}
%\vspace{-2.65in}
%\hspace{-.7in}
%\includegraphics[scale=.75]{ppaf3pb}
%\vspace{-2.8in}
\begin{center}
\includegraphics[scale=.7]{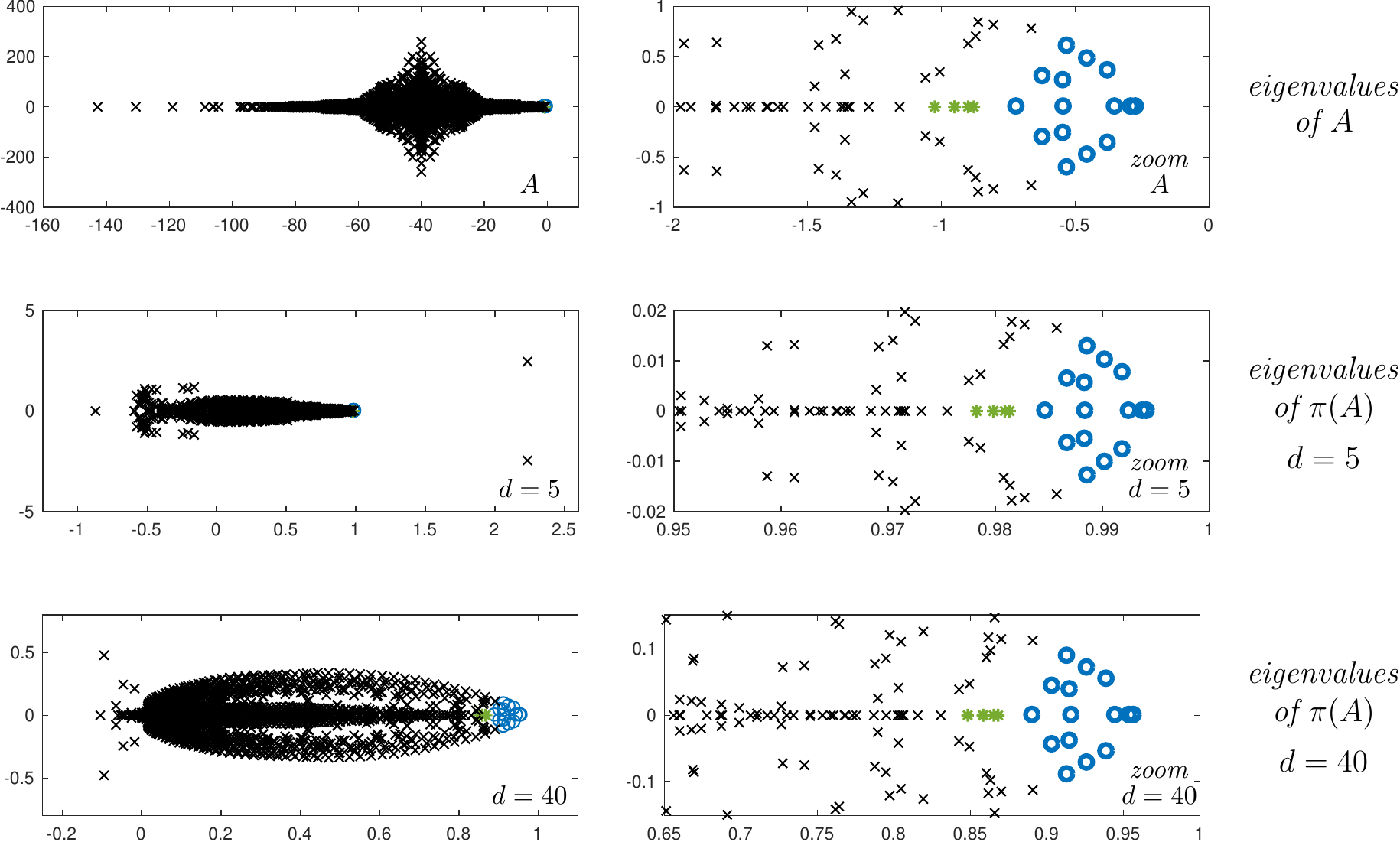}
\end{center}
\caption{\label{fig:ppaf3pb}
Example~2 (matrix Af23560).  The top plots show the eigenvalues of $A$, with the $\nev=15$ desired eigenvalues (nearest~0) as blue circles and the $k-\nev=5$ buffer eigenvalues as green stars; the undesired eigenvalues are shown as black pluses.  The middle plots show the eigenvalues of the preconditioned matrix $\pi(A)$ for degree $\deg=5$; the desired eigenvalues are mapped near~1, but these are interior eigenvalues of $\pi(A)$.
The bottom plots show the eigenvalues of $\pi(A)$ for $\deg=40$: the desired eigenvalues are now on the exterior of the spectrum.}
\end{figure}

%%%%%%%%%%%%%%%%%%%%%%%%%%%%%%%%%%%%%%%%%%%%%%%%%%%%%%%%%%%%%%%%%%%%%%%%%%%%%%%%
{\it Example 2.}
Now consider the matrix Af23560 from Matrix Market, a nonsymmetric matrix of size $n = \mbox{23,560}$ with an average of $\nnzr\approx 19.6$~nonzeros per row.  Again, we seek the smallest magnitude eigenvalues.  The spectrum is complex; see the top left of Figure~\ref{fig:ppaf3pb}.  Polynomial preconditioning is less effective for this example than the previous one, since low degree polynomials struggle to sufficiently isolate the smallest magnitude eigenvalues.  With less sparsity than the last example, the matrix-vector products form a larger part of the expense for this example.

\begin{figure}
%\vspace{-2.6in}
%\hspace{-.7in}
%\includegraphics[scale=.75]{ppaf3}
%\vspace{-2.8in}
\begin{center}
\includegraphics[scale=.63]{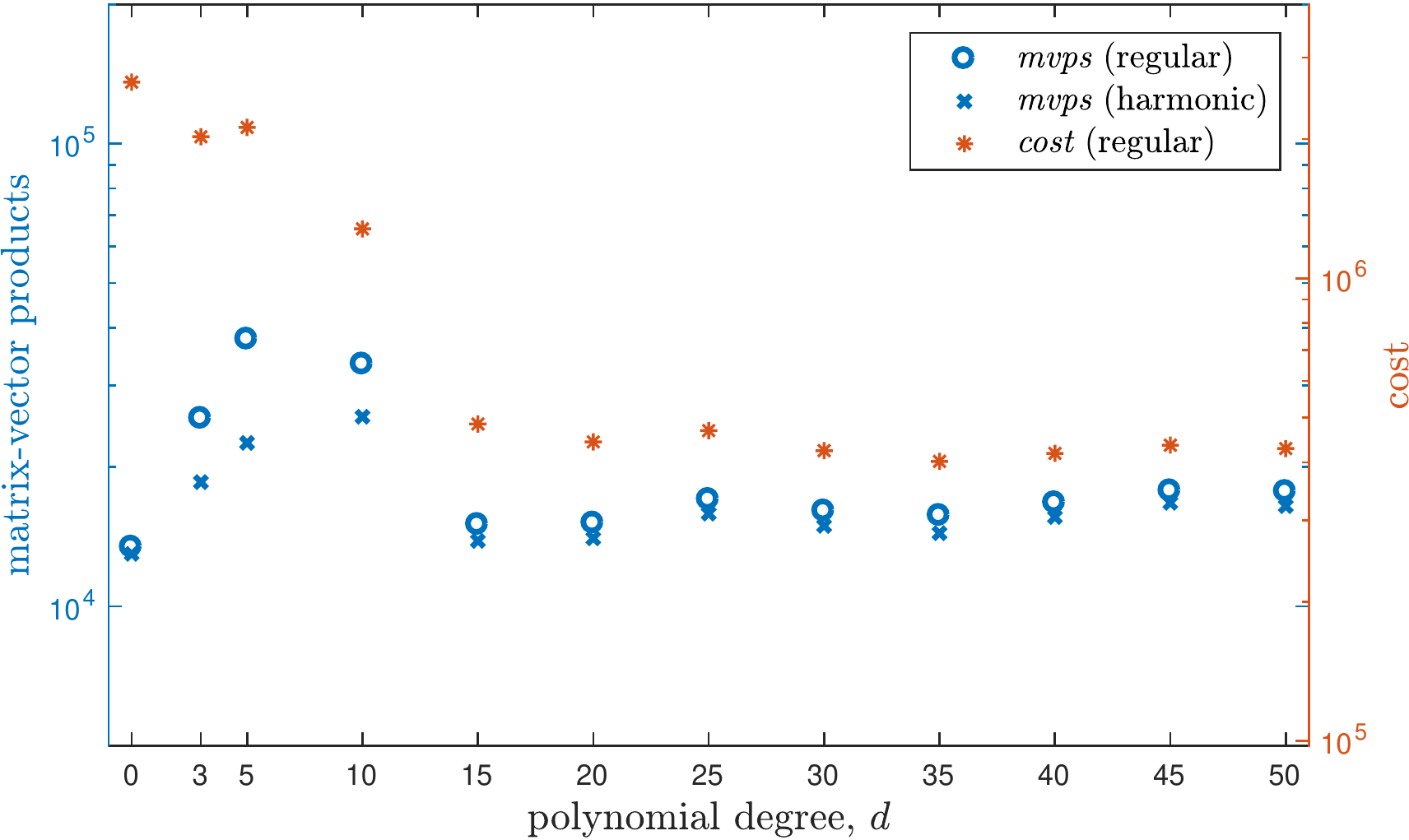}
\end{center}
\vspace*{-5pt}
\caption{\label{fig:ppaf3}
Example~2 (matrix Af23560):
A comparison of convergence between standard Arnoldi(50,20) ($\deg=0$) and polynomial preconditioned Arnoldi(50,20) with $\deg = 3, 5, 10, 15, \ldots, 50$.
Circles indicate the matrix-vector products for regular Arnoldi restarted with Ritz vectors; crosses show the matrix-vector products for Arnoldi restarted with \emph{harmonic} Ritz vectors (left axis).  The approximate cost ($\cost=\nnzr\times\mvps+\vops$) for regular Arnoldi is shown by asterisks (right axis).}
\end{figure}

Figure~\ref{fig:ppaf3} compares convergence for different degree polynomials to no preconditioning ($\deg=0$), showing both matrix-vector products (left axis) and the cost estimate $\cost = \nnzr \times \mvps + \vops$ (right axis).  The right axis  is scaled by~$\nnzr\approx 19.6$, relative to the left axis, to make it easy to see the proportion of cost due to matrix-vector products.  
Figure~\ref{fig:ppaf3} also shows the matrix-vector products required by the harmonic Arnoldi method~\cite{IE,HRAM,PaPavdV} (which keeps harmonic Ritz vectors when restarting, rather than standard Ritz vectors).  (The associated cost estimate is not shown, but is brought down by a similar proportion.)  
Each point in Figure~\ref{fig:ppaf3} is the average of 10~trials.  We note that while all trials converged to the desired tolerance, some only found a subset of the $\nev=15$ desired (smallest magnitude eigenvalues); a few trials found as few as 11~of these desired eigenvalues.

\begin{figure}
%\vspace{-2.65in}
%\hspace{-.7in}
%\includegraphics[scale=.75]{ppaf3pa}
%\vspace{-2.8in}
\begin{center}
\includegraphics[scale=.7]{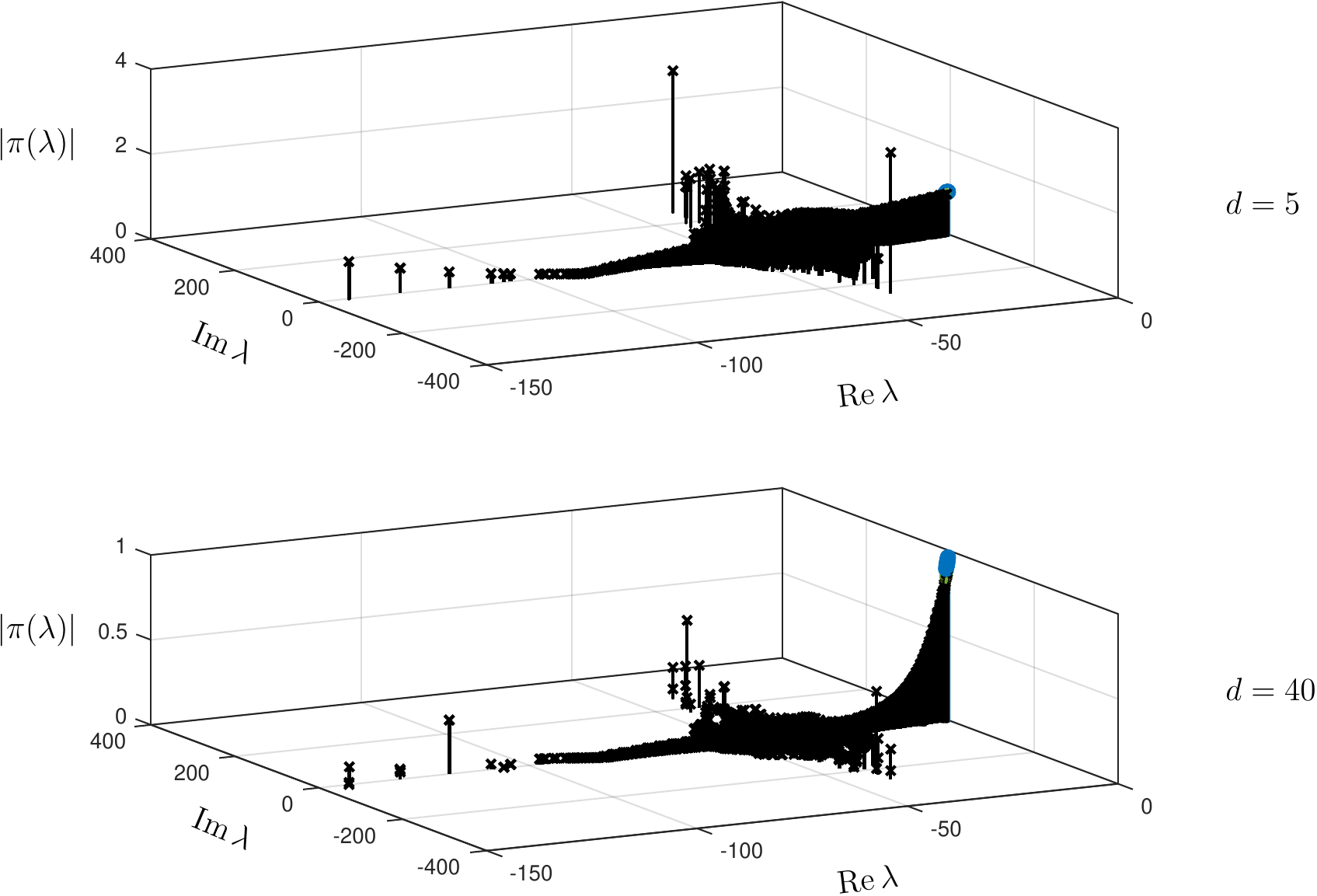}
\end{center}
\caption{\label{fig:ppaf3pa}
Example 2  (matrix  Af23560):  The magnitude of the GMRES polynomial $\pi$ at each eigenvalue of $A$, for degrees $\deg=5$ (top) and $\deg=40$ (bottom).}
\end{figure}

Focusing on $\deg = 5$ (with standard, not harmonic, Arnoldi), polynomial preconditioning uses 2.81~times as many matrix-vector products as with no preconditioning.  Nevertheless, the cost estimate is reduced by about 20\%, because the preconditioned method uses many fewer orthogonalization steps.  For $\deg = 40$, the matrix-vector product count is nearly 25\% higher than for no preconditioning, but $\cost$ is reduced by 84\%.  
Figure~\ref{fig:ppaf3pa} shows the magnitude of the $\deg=5$ and $\deg=40$ minimum residual polynomials $\mrpoly(\lambda)$ on the eigenvalues $\lambda$ of $A$ in the complex plane for a typical example.  In the top plot ($\deg=5$) this polynomial is~1 at the origin, $\mrpoly(0)=1$, but the degree is not large enough for $|\mrpoly(\lambda)|$ to be small at all the eigenvalues of $A$: in particular, $|\mrpoly(\lambda)|$ is large at some large-magnitude eigenvalues on the periphery of the spectrum, with $|\mrpoly(\lambda)|>3$ for one conjugate pair of eigenvalues.
The middle left portion of Figure~\ref{fig:ppaf3pb} shows the resulting spectrum of $\mrpoly(A)$.  
We seek the eigenvalues of $A$ near the origin, and expect $\mrpoly$ to map these close to~1; this is the case, but now these eigenvalues are actually in the interior of the spectrum, due to other eigenvalues $\lambda$ for which $|\mrpoly(\lambda)|>1$, complicating the eigenvalue computation.
This situation motivates our use of the harmonic Arnoldi method (with the same $\deg=5$ polynomial preconditioning), since the harmonic Arnoldi variant is often better for interior eigenvalue problems~\cite{HRAM}.  
For the experiments shown in Figure~\ref{fig:ppaf3}, the harmonic Arnoldi method reduces the matrix-vector products by 40\%.
The lowest portions of Figures~\ref{fig:ppaf3pa} and~\ref{fig:ppaf3pb} show similar spectral information for the minimum residual polynomial with $\deg = 40$.
In this case $\mrpoly$ does a better job of being small at the eigenvalues not near the origin, but it still has some trouble at the eigenvalues with largest imaginary part, which get mapped close to~1.
 Nevertheless, these eigenvalues are no longer mapped so that an interior eigenvalue problem is created.  The right side of Figure~\ref{fig:ppaf3pb} shows close-up views of the  spectrum of $A$ on top and the polynomial preconditioned spectra of $\mrpoly(A)$ for $\deg=5$ and $\deg=40$ below.  Though the overall spectra are vastly changed by the polynomial preconditioning, these close-ups are very similar.

%%%%%%%%%%%%%%%%%%%%%%%%%%%%%%%%%%%%%%%%%%%%%%%%%%%%%%%%%%%%%%%%%%%%%%%%%%%%%%%%
%%%%%%%%%%%%%%%%%%%%%%%%%%%%%%%%%%%%%%%%%%%%%%%%%%%%%%%%%%%%%%%%%%%%%%%%%%%%%%%%

%%%%%%%%%%%%%%%%%%%%%%%%%%%%%%%%%%%%%%%%%%%%%%%%%%%%%%%%%%%%%%%%%%%%%%%%%%%%%%%%
%%%%%%%%%%%%%%%%%%%%%%%%%%%%%%%%%%%%%%%%%%%%%%%%%%%%%%%%%%%%%%%%%%%%%%%%%%%%%%%%

%%%%%%%%%%%%%%%%%%%%%%%%%%%%%%%%%%%%%%%%%%%%%%%%%%%%%%%%%%%%%%%%%%%%%%%%%%%%%%%%
%%%%%%%%%%%%%%%%%%%%%%%%%%%%%%%%%%%%%%%%%%%%%%%%%%%%%%%%%%%%%%%%%%%%%%%%%%%%%%%%

\begin{figure}
%\vspace{-2.65in}
%\hspace{-.7in}
%\includegraphics[scale=.75]{ppaf6pb}
%\vspace{-2.8in}
%\caption{\label{fig:ppaf6pb}
%The matrix E20r0100.  Eigenvalues are shown for the original matrix (top plots) and polynomial preconditioned matrices with degree $\deg=10$ (middle plots) and degree $\deg=25$ (bottom plots).  Closeups of the desired eigenvalues are shown on the right side of the figure, with circles for the~15 nearest the desired value (0~or the original spectrum and 1~for the preconditioned spectra), asterisks for the next~5 and pluses for the others.}
\begin{center}
\includegraphics[scale=.7]{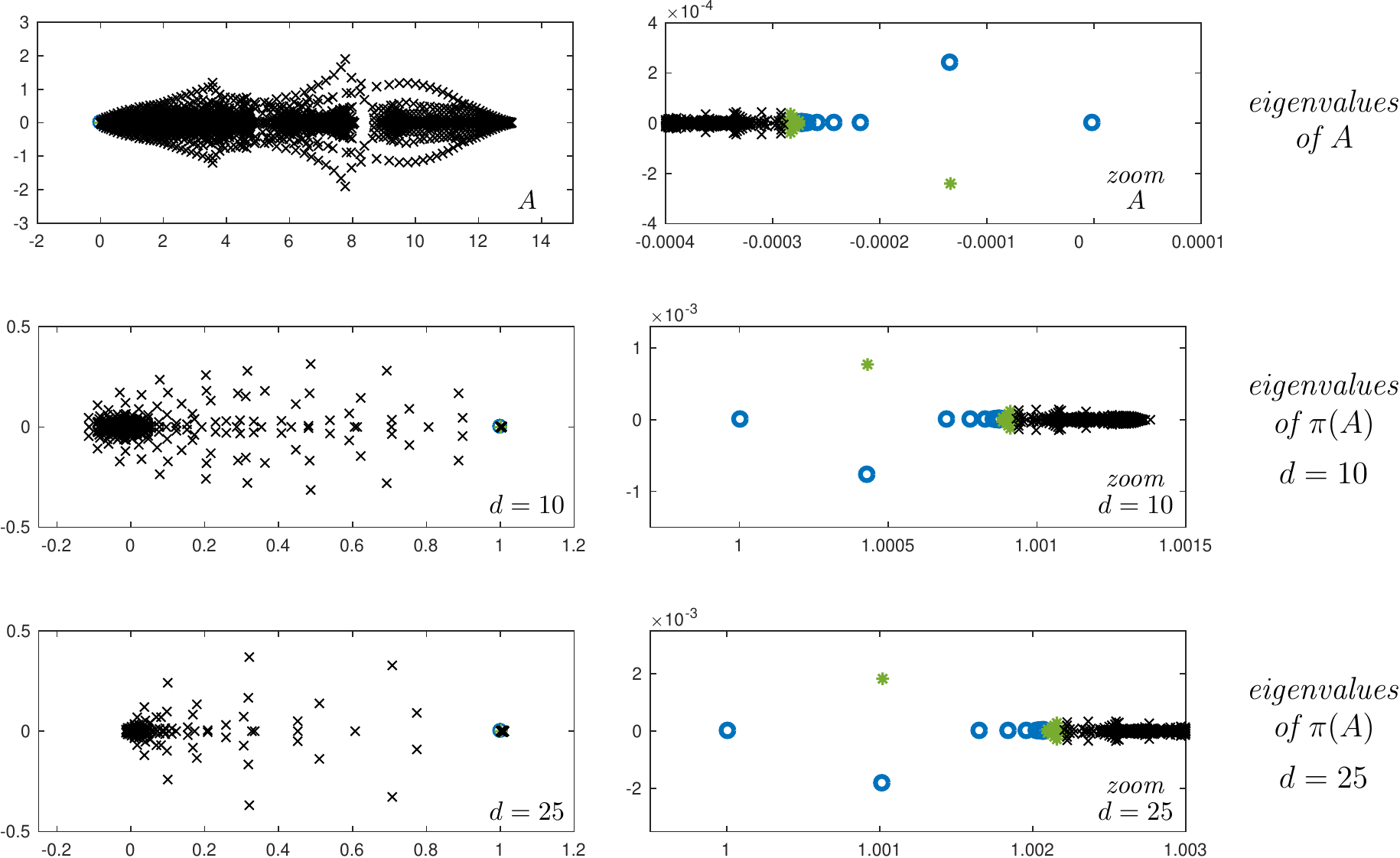}
\end{center}
\caption{\label{fig:ppaf6pb}
Example~3 (matrix E20r0100).  The top plots show the eigenvalues of $A$, with the $\nev=15$ desired eigenvalues (nearest~0) as blue circles and the $k-\nev=15$ buffer eigenvalues as green stars; the undesired eigenvalues are shown as black pluses.  The middle plots show the eigenvalues of the preconditioned matrix $\pi(A)$ for degree $\deg=10$; the desired eigenvalues are mapped near~1. 
The bottom plots show the eigenvalues of $\pi(A)$ for $\deg=25$.}
\end{figure}

{\it Example 3.} The matrix E20r0100 from Matrix Market also has a complex spectrum, but this time polynomial preconditioning is more effective, reducing the overall matrix-vector products.  This matrix has dimension $n=4241$ with an average of nearly 31~nonzeros per row.  As before, we seek the~$\nev=15$ eigenvalues nearest the origin, which  are in the interior of the spectrum ($A$ has 1199~eigenvalues with (quite small) negative real parts), though now we use Arnoldi(100,30) iterations to access larger subspaces.  Figure~\ref{fig:ppaf6pb} shows how polynomial preconditioning changes the spectrum, and Table~\ref{fig:ppaf6pb} gives results for different degree polynomials.  A degree $\deg=15$ polynomial reduces the number of matrix-vector products by a factor of almost~9 and vector operations by a factor of more than~125.  Matrix-vector products are not further reduced with higher degree polynomials, but the vector operations are.   

%%%%%%%%%%%%%%%%%%%%%%%%%%%%%%%%%%%%%%%%%%%%%%%%%%%%%%%%%%%%%%%%%%%%%%%%%%%%%%%%
%%%%%%%%%%%%%%%%%%%%%%%%%%%%%%%%%%%%%%%%%%%%%%%%%%%%%%%%%%%%%%%%%%%%%%%%%%%%%%%%

%%%%%%%%%%%%%%%%%%%%%%%%%%%%%%%%%%%%%%%%%%%%%%%%%%%%%%%%%%%%%%%%%%%%%%%%%%%%%%%%
\begin{table}
\caption{\label{tbl:ppaf6}
Example~3 (Matrix e20r0100, $n=4241$):
results for polynomial preconditioning with different degree polynomials using Arnoldi(100,30) to seek $\nev=15$ eigenvalues, averaged over 10~trials.}

\vspace*{-5pt}
\begin{center}
\begin{tabular}{|r|r|r|r|r|r|r|} \hline\hline 
\multicolumn{1}{|c|}{degree} & \multicolumn{1}{|c|}{cycles} & \multicolumn{1}{|c|}{$\mvps$} & \multicolumn{1}{|c|}{$\vops$} & \multicolumn{1}{|c|}{dot products} & \multicolumn{1}{|c|}{$\cost$} & \multicolumn{1}{|c|}{\# correct} \\ 
\multicolumn{1}{|c|}{$\deg$} & \multicolumn{1}{|c|}{ } & \multicolumn{1}{|c|}{} & \multicolumn{1}{|c|}{(millions)} & \multicolumn{1}{|c|}{(millions)} & \multicolumn{1}{|c|}{(millions)} & \multicolumn{1}{|c|}{eigenvalues} \\ \hline 
  0\ \ \ &  \mbox{7,959.1} &  \mbox{558,258.0} &    171.80\ \ \ &     74.41 \ \ \ \ \ &    189.10\ \ \ &  14.8\ \ \ \ \ \\ \hline 
  5\ \ \ &   308.6 &  \mbox{109,067.6} &      6.82\ \ \ &      2.89 \ \ \ \ \ &     10.20\ \ \ &  14.0\ \ \ \ \ \\ \hline 
 10\ \ \ &    97.4 &   \mbox{68,926.6} &      2.19\ \ \ &      0.91 \ \ \ \ \ &      4.33\ \ \ &  14.0\ \ \ \ \ \\ \hline 
 15\ \ \ &    58.9 &   \mbox{62,740.8} &      1.35\ \ \ &      0.55 \ \ \ \ \ &      3.29\ \ \ &  14.0\ \ \ \ \ \\ \hline 
 20\ \ \ &    46.0 &   \mbox{65,510.1} &      1.07\ \ \ &      0.43 \ \ \ \ \ &      3.10\ \ \ &  14.0\ \ \ \ \ \\ \hline 
 25\ \ \ &    43.2 &   \mbox{76,936.4} &      1.02\ \ \ &      0.41 \ \ \ \ \ &      3.41\ \ \ &  14.0\ \ \ \ \ \\ \hline 
 50\ \ \ &    44.5 &  \mbox{158,253.2} &      1.14\ \ \ &      0.42 \ \ \ \ \ &      6.04\ \ \ &  14.7\ \ \ \ \ \\ \hline 
 75\ \ \ &    27.9 &  \mbox{149,730.0} &      0.77\ \ \ &      0.27 \ \ \ \ \ &      5.41\ \ \ &  14.8\ \ \ \ \ \\ \hline 
100\ \ \ &    14.8 &  \mbox{107,274.7} &      0.46\ \ \ &      0.15 \ \ \ \ \ &      3.78\ \ \ &  14.0\ \ \ \ \ \\ \hline 
\end{tabular} 
\end{center}
\end{table}

\section{Two Starting Vectors} \label{sec:twob}
%%%%%%%%%%%%%%%%%%%%%%%%%%%%%%%%%%%%%%%%%%%%%%%%%%%%%%%%%%%%%%%%%%%%%%%%%%%%%%%%

For the examples in the last section we randomly generated the starting vector used to create $\mrpoly$, an approach that seems to work quite well in practice.  However, here we consider the possibility that an unusual or skewed starting vector for $\mrpoly$ can give bad results, and show how two starting vectors can be used to generate $\pi$ to minimize the risk of a bad starting vector.

{\it Example 4.}  Let $A$ be  diagonal  with main diagonal elements $1, 2, 3, \ldots, 1000$.  We run Arnoldi(50,20) to calculate $\nev=15$ eigenvalues to residual norm of $10^{-8}$ and use a different starting vector for the polynomial preconditioned Arnoldi loop than was used for developing the polynomial (with $\deg=10$).  
It takes only one cycle to find all 15~correct eigenpairs.  Next, we make the last 100~components of the starting vector for the polynomial small by multiplying them by~0.01. 
The resulting $\mrpoly$ is not small much past $\lambda=930$;  $\mrpoly(\lambda)$ goes up to at least~7 at $\lambda=1000$, transforming the problem of finding the eigenvalues near~1 into an interior eigenvalue problem.  Convergence is much slower, with 16.8~cycles needed (average of 10~trials), and then only the first four desired eigenvalues ($\lambda=1,2,3,4$) are found.  The remaining computed eigenvalues fall in $\{938,939,\ldots,956\}$, depending on the run:  eigenvalues $\lambda$ that $\mrpoly$ maps closer to the target point~1 than the desired eigenvalues $\lambda=5,\ldots,15$.

We give an algorithm that uses two starting vectors to determine one polynomial, applying GMRES to a $2\times2$ block diagonal system of dimension $2n$.

\vspace{.07in}
\begin{center}
\textbf{Polynomial Determined by Two Starting Vectors}
\end{center}
\begin{enumerate}
 \item {\bf Set-up:} Generate random vectors $b_1$ and $b_2$, with $\|b_1\|=\|b_2\|=1/\sqrt{2}$.  Let 
 \[ \widehat{b} = \begin{bmatrix}
b_1 \\
b_2
\end{bmatrix}, \qquad
\widehat{A} = \begin{bmatrix}
A & 0 \\
0 & A
\end{bmatrix}.\]
 \item {\bf Generate polynomial:} Run a cycle of GMRES($\deg$) with starting vector $\widehat{b}$ and matrix $\widehat{A}$, and find the roots of the GMRES polynomial $\mrpoly$.
\end{enumerate}
\vspace{.08in}

This approach essentially uses two Krylov subspaces, one each with $b_1$ and $b_2$, and so takes into account both starting vectors.  We tested Example~4 with the skewed starting vector for $b_1$, but $b_2$ a random vector.  The results are good: we now need only one or two cycles to compute all the desired small eigenvalues.

%%%%%%%%%%%%%%%%%%%%%%%%%%%%%%%%%%%%%%%%%%%%%%%%%%%%%%%%%%%%%%%%%%%%%%%%%%%%%%%%%%%%%%%%%%%%%%%%%%
\section{Damped Polynomials}
%%%%%%%%%%%%%%%%%%%%%%%%%%%%%%%%%%%%%%%%%%%%%%%%%%%%%%%%%%%%%%%%%%%%%%%%%%%%%%%%%%%%%%%%%%%%%%%%%%

Sometimes the GMRES polynomial $\pi$ is not an ideal preconditioner for eigenvalue calculations.
In this section we examine several scenarios that can lead to poor preconditioners, and propose techniques to tame the extreme behavior of the GMRES polynomial; we call the result a ``damped'' polynomial.  Section~\ref{sec:overenthused} addresses the case were $\mrpoly(\lambda)$ is too large at undesired eigenvalues $\lambda$; Section~\ref{sec:tooeasy} investigates problems that are \emph{too easy}, making $\mrpoly(\lambda)$  too small at desired $\lambda$.  
Both cases can be improved with a damped polynomial.

%%%%%%%%%%%%%%%%%%%%%%%%%%%%%%%%%%%%%%%%%%%%%%%%%%%%%%%%%%%%%%%%%%%%%%%%%%%%%%%%%%%%%%%%%%%%%%%%%%
\subsection{Overenthusiastic polynomials} \label{sec:overenthused}
%%%%%%%%%%%%%%%%%%%%%%%%%%%%%%%%%%%%%%%%%%%%%%%%%%%%%%%%%%%%%%%%%%%%%%%%%%%%%%%%%%%%%%%%%%%%%%%%%%

The next example considers a case where $|\mrpoly(\lambda)|$ is large at unwanted eigenvalues; it also shows how polynomial preconditioning can be effective even for a matrix with less sparsity than the earlier examples.

\begin{figure}
%\vspace{-2.65in}
%\hspace{-.7in}
%\includegraphics[scale=.75]{ppaf5.pdf}
%\vspace{-2.8in}
\begin{center}
\includegraphics[scale=.63]{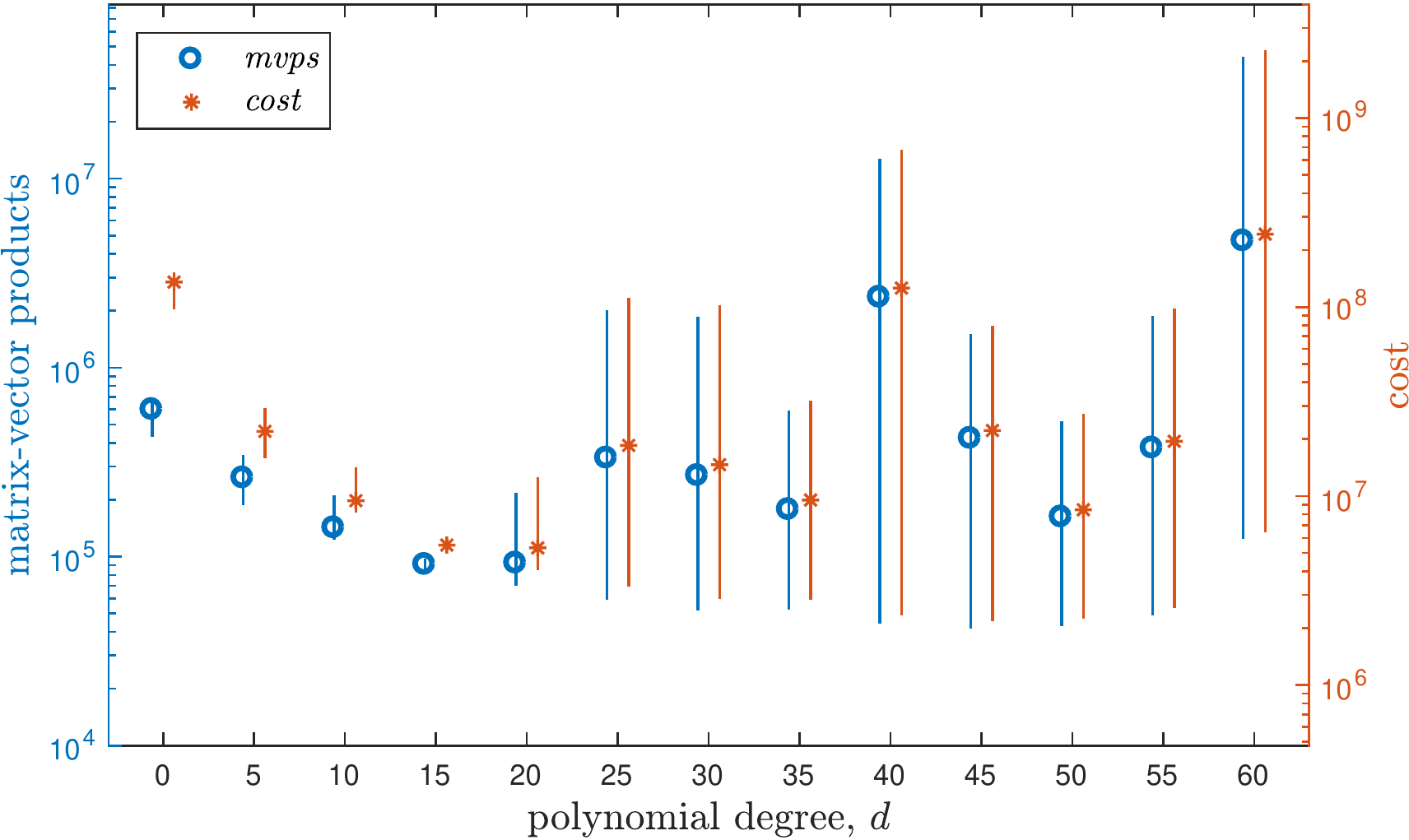}
\end{center}

\vspace*{-5pt}
\caption{\label{fig:ex5_dplot}
Example~5 (S1rmq4m1 matrix): Convergence of standard Arnoldi(50,20) ($d=0$) and polynomial preconditioned Arnoldi(50,20) with $d=5,10,15,\ldots, 60$.  Blue circles indicate the average number of matrix-vector products (left axis); red asterisks show the average approximate cost ($\cost = \nnzr\times \mvps+\vops$, right axis).  These results vary widely over ten trials;  vertical bars show the minimum and maximum matrix-vector products and cost over these trials.}
\end{figure}

{\it Example 5.} Consider S1rmq4m1 from Matrix Market, a symmetric, positive definite matrix of size $n = 5489$ with an average of $\nnzr\approx 47.8$ nonzeros per row.  Finding the small eigenvalues is difficult because they are packed closely together relative to the whole spectrum. The first 15 range from $0.3797$ to $1.9027$, the 21st is 2.2630 and the largest is $6.87 \times 10^{5}$.  Figure~\ref{fig:ex5_dplot} shows the performance of Arnoldi(50,20) applied to find the $\nev=15$ smallest magnitude eigenvalues of $A$.  
The gap ratio $\frac{|\lambda_{21}-\lambda_{15}|}{|\lambda_{5489}-\lambda_{21}|}\approx 5.24\times 10^{-7}$ roughly describes the eventual convergence for the 15th eigenvalue.
For a typical case of polynomial preconditioning with $\deg=20$, 
$\frac{|\mrpoly(\lambda_{21})-\mrpoly(\lambda_{15})|}{|\mrpoly(\lambda_{2737})-\mrpoly(\lambda_{21})|} = \frac{|0.99728 - 0.99771|}{|-1.82056 - 0.99728|} \approx 1.54\times 10^{-4}$ (though recall each preconditioned Arnoldi iteration requires 20~matrix-vector products).
Nevertheless, there is an improvement in matrix-vector products by a factor of~6.5 over these ten trials, and $\cost = \nnzr\times\mvps + \vops$ is reduced by a factor of~25.51.
Though $A$ has a relatively large number of nonzeros, $\nnzr\approx 47.8$, the cost of vector operations still dominates for regular Arnoldi;
in contrast, with polynomial preconditioning the matrix-vector products are the bigger expense.
Moreover, for $\deg=20$ polynomial preconditioning decreases the dot products
by a factor of~130.34 over regular Arnoldi.
%; when $\deg=50$, by a factor of \hot{over 800}.  % skip, given variability in large d for this example

\begin{comment}
                       lambda_1 :       0.37969627 
                       lambda_15:       1.90275206 
                       lambda_21:       2.26295090 
                        lambda_n:   6.87432150e+05 
               regular gap ratio:      5.24e-07 

 preconditioned gap ratio, d=20:      1.54e-04 
  index for d=20 gap ratio (15):    15
  index for d=20 gap ratio (21):    21
  index for d=20 gap ratio (n):   2737

pi(lam) for d=20 gap ratio (15):       0.99771
pi(lam) for d=20 gap ratio (21):       0.99728
pi(lam) for d=20 gap ratio (n):       -1.82056

 matrix vector products, regular:     597494.00 
 matrix vector products, d=20:         91966.50 
    factor of matvec improvement:          6.50 

          cost estimate, regular:  135598309.19 
          cost estimate, d=20:       5315906.17 
      factor of cost improvement:         25.51 

          dots estimate, regular:   44057540.10 
          dots estimate, d=20:        338010.20 
      factor of dots improvement:        130.34 
\end{comment}

%%%%%%%%%%%%%%%%%%%%%%%%%%%%%%%%%%%%%%%%%%%%%%%%%%%%%%%%%%%%%%%%%%%%%%%%%%%%%%%%%%%%%%%%%%%%%%%%%%
%%%%%%%%%%%%%%%%%%%%%%%%%%%%%%%%%%%%%%%%%%%%%%%%%%%%%%%%%%%%%%%%%%%%%%%%%%%%%%%%%%%%%%%%%%%%%%%%%%

\begin{figure}
%\vspace{-2.7in}
%\hspace{-.7in}
%\includegraphics[scale=.75]{ppaf5p.pdf}
%\vspace{-2.8in}
\begin{center}
\includegraphics[width=3.75in]{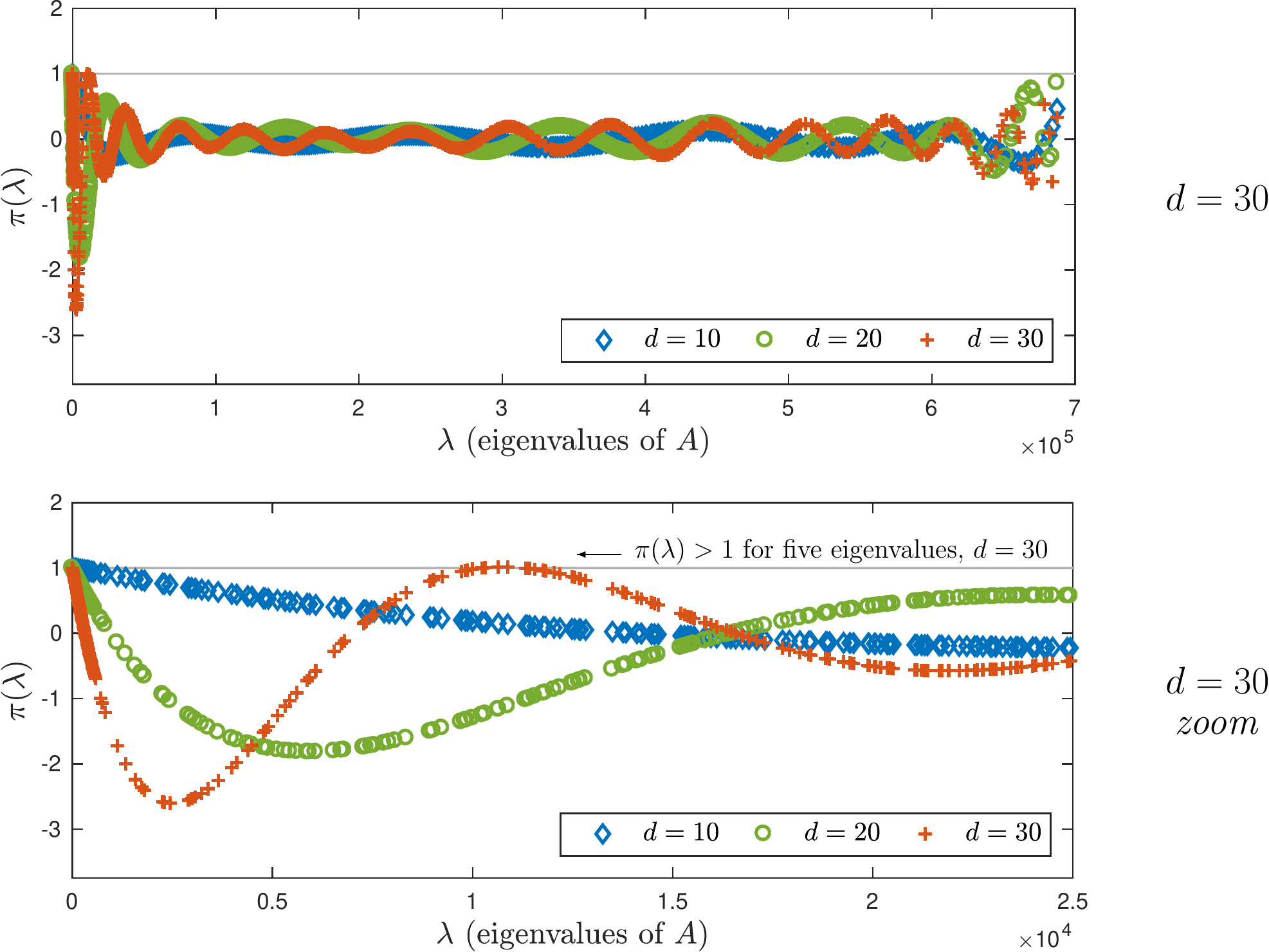}
\end{center}
\vspace*{-4pt}
\caption{\label{fig:ex5_poly1}
Example~5 (S1rmq4m1 matrix):  Polynomials of degree $\deg=10$, 20 and 30.  The top plot shows $\mrpoly(\lambda)$ for all eigenvalues $\lambda$ of $A$; the bottom plot zooms in near the origin.}
\end{figure}

Figure~\ref{fig:ex5_dplot} shows that polynomials of degree $d > 20$ can cause problems. 
The performance starts to vary widely over our 10~trial runs; for each $d>20$, at least one run failed to find the correct $\nev=15$ smallest eigenvalues.  Figure~\ref{fig:ex5_poly1} indicates the problem, showing $\pi(\lambda)$ for $\deg=10$, 20, and 30 for the first of our ten trials.  The slope of $\pi$ at the origin is much steeper for $\deg=20$ than for $\deg=10$, explaining the faster convergence.  Degree $\deg=30$ is even steeper, but has a problem near $\mbox{10,000}$, where the $\mrpoly(\lambda)>1$ for five eigenvalues:  $\pi$ maps eigenvalues of $A$ from the interior of the spectrum to the exterior of $\mrpoly(A)$, mixing seven of them amongst or above the $\nev=15$ desired smallest eigenvalues of $\mrpoly(A)$ near~1.  The resulting interior eigenvalue problem can lead to slow convergence and spurious eigenvalues.  Of the converged Ritz values for this example, only~11 fall among the $\nev$ desired smallest magnitude eigenvalues; the others are unwanted interior eigenvalues of $A$.   As Figure~\ref{fig:ex5_dplot} shows, erratic convergence continues for larger values of $\deg$.   (One run with $\deg=60$ failed to compute \emph{any} of the eigenvalues correctly, and took 24,287~cycles to converge; in a  run not included in the plot, a $\deg=30$ run failed to converge in 30,000~cycles.) 

We call the polynomials that jump up too high ``overenthusiastic".
This matrix seems prone to such polynomials because about half its spectrum is near 0 (2703 eigenvalues are less than 600; the next 2786 go from 700 up to $6.87 \times 10^{5}$): the GMRES polynomial seems to concentrate on these small eigenvalues, while focusing less precisely on the others.

%%%%%%%%%%%%%%%%%%%%%%%%%%%%%%%%%%%%%%%%%%%%%%%%%%%%%%%%%%%%%%%%%%%%%%%%%%%%%%%%%%%%%%%%%%%%%%%%%%
%%%%%%%%%%%%%%%%%%%%%%%%%%%%%%%%%%%%%%%%%%%%%%%%%%%%%%%%%%%%%%%%%%%%%%%%%%%%%%%%%%%%%%%%%%%%%%%%%%

\begin{figure}
%\vspace{-2.7in}
%\hspace{-.7in}
%\includegraphics[scale=.75]{ppaf5p2.pdf}
%\vspace{-2.8in}
\begin{center}
\includegraphics[width=3.75in]{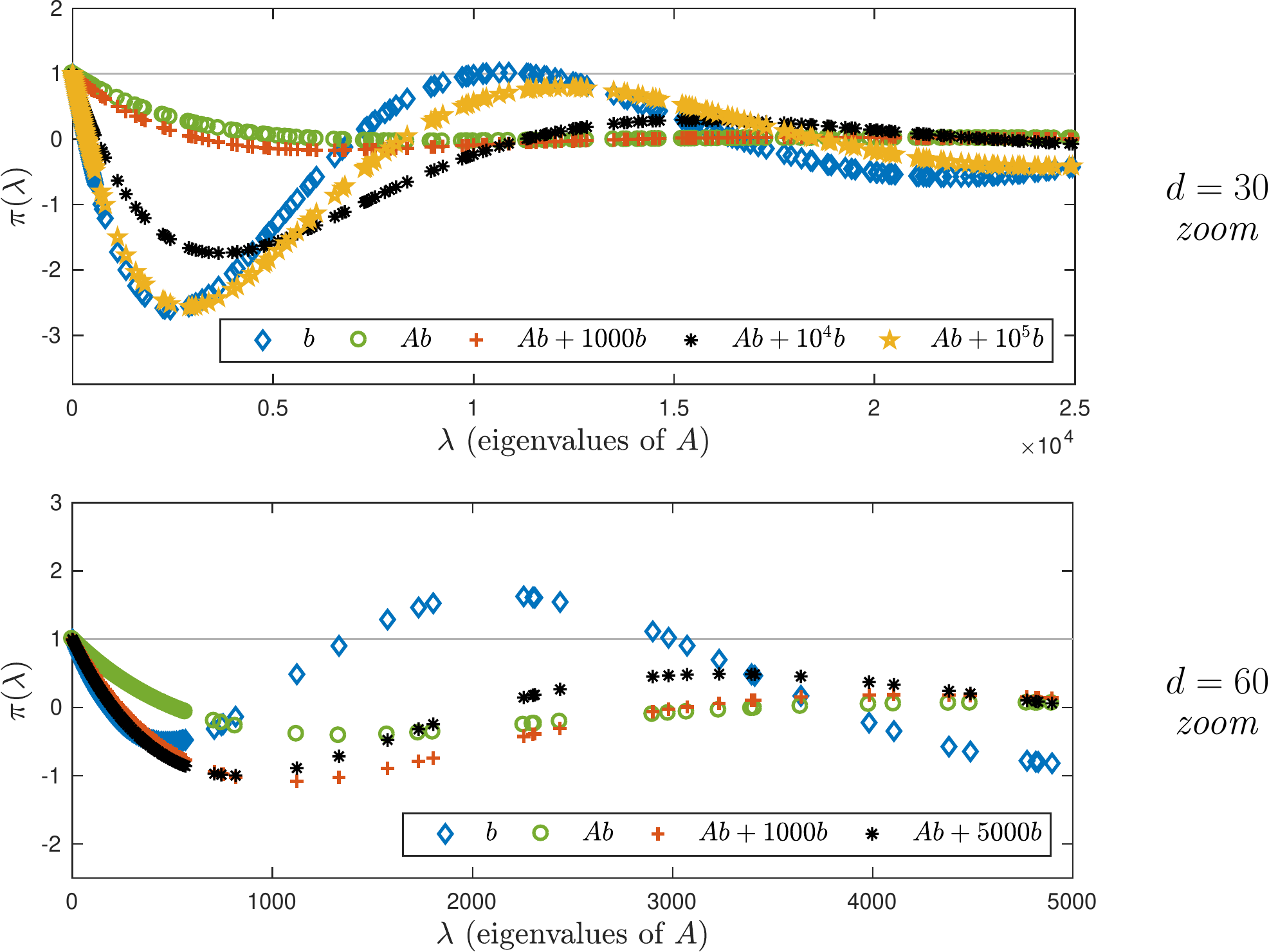}
\end{center}
\vspace*{-4pt}
\caption{\label{fig:ex5_poly2}
Example 5 (S1rmq4m1 matrix): Close-ups near the origin of polynomials of degree $\deg=30$ (top) and $\deg=60$ (bottom), using different damping strategies.}
\vspace*{-12pt}
\end{figure}

Damping the polynomial provides a possible remedy to overenthusiasm.  
Damping techniques are considered in~\cite{LiXiVeYaSa} for symmetric matrices and for polynomials that approximate the Dirac delta function.
Here we take a different approach.  We first damp by changing the starting vector for GMRES($\deg$) from a random vector $b$ to $A@b$: premultiplication by $A$ will generally reduce the components of $b$ in the eigenvectors corresponding to the small eigenvalues.  (Think of performing one step of the power method.)   The GMRES polynomial for $A@b$ is then less likely to be overenthusiastic, because it does not try as hard to be small at the small eigenvalues.  The top portion of Figure~\ref{fig:ex5_poly2} shows polynomials of degree $\deg=30$.  Diamonds show the GMRES polynomial generated from $b$; circles show the damped GMRES polynomial generated from $A@b$.  The damped polynomial no longer jumps too high in the middle of the spectrum, but its slope is much less steep at the origin than for the standard polynomial: thus it yields slower than desired convergence. For the run shown here,  $\cost= 4.28\times10^6$, about the same as for undamped $\deg=20$ ($\cost=4.23\times10^6$); of course, both are better than $\cost=97.4\times 10^6$ required without polynomial preconditioning.
Next we try starting vectors of the form $A@b + \alpha b$.  Figure~\ref{fig:ex5_poly2}
shows results for a few choices of $\alpha$.  For $\alpha = 10^5$,
$\cost = 3.34\times 10^{6}$, but this $\mrpoly$ is on the verge of being overenthusiastic.

The bottom of Figure~\ref{fig:ex5_poly2} repeats this experiment with  $\deg=60$.  The undamped polynomial is very enthusiastic; it goes far above~1.0.  Again this can be controlled with the starting vector $A@b + \alpha b$, but now a higher proportion of the $Ab$ term is needed. 
With $\alpha = 5000$, the $\cost$ goes down to $2.50\times10^6$.

%%%%%%%%%%%%%%%%%%%%%%%%%%%%%%%%%%%%%%%%%%%%%%%%%%%%%%%%%%%%%%%%%%%%%%%%%%%%%%%%%%%%%%%%%%%%%%%%%%
%%%%%%%%%%%%%%%%%%%%%%%%%%%%%%%%%%%%%%%%%%%%%%%%%%%%%%%%%%%%%%%%%%%%%%%%%%%%%%%%%%%%%%%%%%%%%%%%%%

%%%%%%%%%%%%%%%%%%%%%%%%%%%%%%%%%%%%%%%%%%%%%%%%%%%%%%%%%%%%%%%%%%%%%%%%%%%%%%%%%%%%%%%%%%%%%%%%%%
\subsection{Easy problems} \label{sec:tooeasy}
%%%%%%%%%%%%%%%%%%%%%%%%%%%%%%%%%%%%%%%%%%%%%%%%%%%%%%%%%%%%%%%%%%%%%%%%%%%%%%%%%%%%%%%%%%%%%%%%%%

The GMRES polynomial $\mrpoly$ starts at $\mrpoly(0)=1$ at the origin and tries to drop down quickly to be near zero over the spectrum.
For easy problems (where the spectrum is well separated from the origin) and high degree polynomials, $\mrpoly$ can drop down too fast, mapping some desired eigenvalues near zero, mixed amongst the undesired eigenvalues.  Figure~\ref{fig:polyplot} shows such behavior: for $\deg = 8$ the first~10 eigenvalues of $\mrpoly(A)$ are separated from the others, but the next few  are mixed with the rest of the spectrum.
The preconditioner $\mrpoly$ will be effective if $\mrpoly(\lambda)\approx 1$ for the desired eigenvalues, and $|\mrpoly(\lambda)|$ is small for the undesired ones.  To promote such behavior one can use small $\deg$, or operate on $A-\mu I$ for $\mu$ near the desired eigenvalues.  Here we show how damping can also help.

\medskip
{\it Example 6.}
Let $A$ be the diagonal matrix of dimension $n=\mbox{10,000}$ with diagonal entries $1, 2, 3, \ldots$, $9999, 10000$.  We  run Arnoldi(50,20) to calculate the 
$\nev=15$~smallest magnitude eigenvalues to residual norm $\rtol=10^{-8}$, averaging over ten trials. 
Without preconditioning, the calculation takes 1,625~$mvps$ and 287,282~$\vops$.  With $\deg = 40$, more matrix-vector products are needed (2,775.6~$\mvps$), but other costs are much lower (15,390.8 $\vops$).  Increasing to $\deg = 50$ changes the results dramatically:  $\mrpoly$ is too small at some desired eigenvalues; all 10~trials miss $\lambda= 13, 14, 15$; some runs miss more.  (The various trials compute unwanted eigenvalues among $\{66,67,\ldots, 72\}$.)  Furthermore, convergence is much slower, with many eigenvalues mapped close together: 12,621.7 $\mvps$ are needed.  For $\deg=50$, the GMRES problem is too easy.

We next try the damped polynomial with $A@b$ as the GMRES starting vector for $\deg=50$.  The performance is better; the correct eigenvalues are found in 2,565.0~$\mvps$, a single cycle.
If we continue increasing $\deg$, even the damped polynomial runs into trouble by coming down too quickly.  When $\deg=80$ most of the trials miss at least one eigenalue, but this can be fixed by damping more, using starting vector $A^2@b$.

%%%%%%%%%%%%%%%%%%%%%%%%%%%%%%%%%%%%%%%%%%%%%%%%%%%%%%%%%%%%%%%%%%%%%%%%%%%%%%%%%%%%%%%%%%%%%%%%%%
\subsection{A heuristic toward automation}
%%%%%%%%%%%%%%%%%%%%%%%%%%%%%%%%%%%%%%%%%%%%%%%%%%%%%%%%%%%%%%%%%%%%%%%%%%%%%%%%%%%%%%%%%%%%%%%%%%

We outline an attempt to automatically determine when to damp, for both too easy and overenthusiastic situations.

\vspace{.07in}
\begin{center}
\textbf{Heuristic for Damping \boldmath $b$ and Adjusting $\deg$}
\end{center}

\begin{enumerate}
 \item Apply one Arnoldi$(m,k)$ cycle to $\pi(A)$ with starting vector $b$.
 \item Compute the Ritz values $\nu_1,\ldots, \nu_m$ for $\pi(A)$ and associated (unit-length) Ritz vectors $y_1,\ldots,y_m$ from $\CK_m(\pi(A),b)$, ordered by increasing distance from~1:
% \vspace*{-3pt}
 \[ |1-\nu_1| \le |1-\nu_2| \le \cdots \le |1-\nu_m|.\]
 %\vspace*{-12pt}
 \item Compute Rayleigh--Ritz eigenvalue estimates \emph{for $A$}: $\mu_j \equiv y_j^*Ay_j^{}$.
 \item Check the \emph{Ideal Order Condition}:
 % \vspace*{-3pt}
        \[ |\mu_1| \le |\mu_2| \le \cdots \le |\mu_{\text{\it nev}}|
             < \min_{\text{\it nev}+1\le j\le k} |\mu_j|.\]
 %\vspace*{-10pt}
 \item If the Ideal Order Condition holds, proceed with the computation;\\
       If the Ideal Order Condition fails:\\
       $\bullet$ Generate a new $\mrpoly$ using starting vector $A@b$.  Repeat steps 1--3.\\
       $\bullet$ If the Ideal Order Condition still fails, replace $d$ with $\lfloor d/2\rfloor$ and try again.\\
       $\bullet$ Continue halving $d$ until the condition holds (or $d=1$).
\end{enumerate}
\vspace{.08in}

{\it Example 6 (continued).}
For the matrix in Example~6 with $d = 50$, a cycle of Arnoldi(50,20) applied to $\mrpoly(A)$ leads to these Rayleigh quotients $\mu_1, \ldots, \mu_k$
for $A$:
\[     1,\ 2,  3,\ \ldots, 12,\ 
   \underline{68.83},\ 
   \underline{67.62},\ 
   \underline{68.72},\ 
   71.74,\ 
   67.26,\ 
   67.86,\ 
   84.61,\ 
   72.80.
 \]
(Continuing the Arnoldi(50,20) algorithm with this $\mrpoly(A)$ only finds 10 of the desired $\nev=15$ eigenvalues.)
The Ideal Order Condition requires these Rayleigh quotients
to increase monotonically in magnitude and be dominated by the last $k-\nev=5$ values; the underlined values violate one or both of these requirements, so the heuristic replaces $b$ with $A@b$ and computes a fresh $\mrpoly$.  This damping is successful: Arnoldi(50,20) converges to all $\nev=15$ desired eigenvalues in just one cycle.

\medskip  
We tested 100~runs each of degrees $\deg=30, 31, \ldots, 60$ (a total of 3100 tests).  The problem is sometimes ``too easy'' starting at $\deg = 39$, with increasing prevalence as $\deg$ increases.  Of the 3100~tests, 1425~failed the first damping test.  (We ran Arnoldi(50,20) anyway: \emph{all~1425 cases failed to find all $\nev=15$ correct eigenvalues}.)  With damping (replacing random $b$ by $A b$), all~1425 cases passed the second damping test; Arnoldi(50,20) proceeded to compute all $\nev=15$ desired eigenvalues.

The damping test is not perfect, however.  Of the 1675~tests that passed the first damping test, about 95\% of the cases successfully computed the $\nev=15$ eigenvalues.  The remaining 82~cases failed to find some of the desired eigenvalues (although all found at least 8~correct eigenvalues).  While this simple damping test is not perfect,  in many cases it signals the need for a modified approach.

\begin{comment}
There are no false positives where the test after the first cycle says to begin again with damping even though it is not necessary.  However, we have been able to rig an example to create a false positive.  There is an interesting case with $d = 38$, because even though the correct eigenvectors are computed, convergence takes seven cycles due to some larger eigenvalues $A$ being mapped by $\pi$ close to the 15th eigenvalue of $\pi(A)$ .
\end{comment}

%%%%%%%%%%%%%%%%%%%%%%%%%%%%%%%%%%%%%%%%%%%%%%%%%%%%%%%%%%%%%%%%%%%%%%%%%%%%%%%%%%%%%%%%%%%%%%%%%%
\section{Stability}
%%%%%%%%%%%%%%%%%%%%%%%%%%%%%%%%%%%%%%%%%%%%%%%%%%%%%%%%%%%%%%%%%%%%%%%%%%%%%%%%%%%%%%%%%%%%%%%%%%

Several earlier examples  hinted that high degree preconditioners can be temperamental.  Here we investigate how such polynomials, with their steep slopes, can lead to a computational instability, and propose a way to cope.

\begin{figure}
\begin{center}
\includegraphics[width=3.25in]{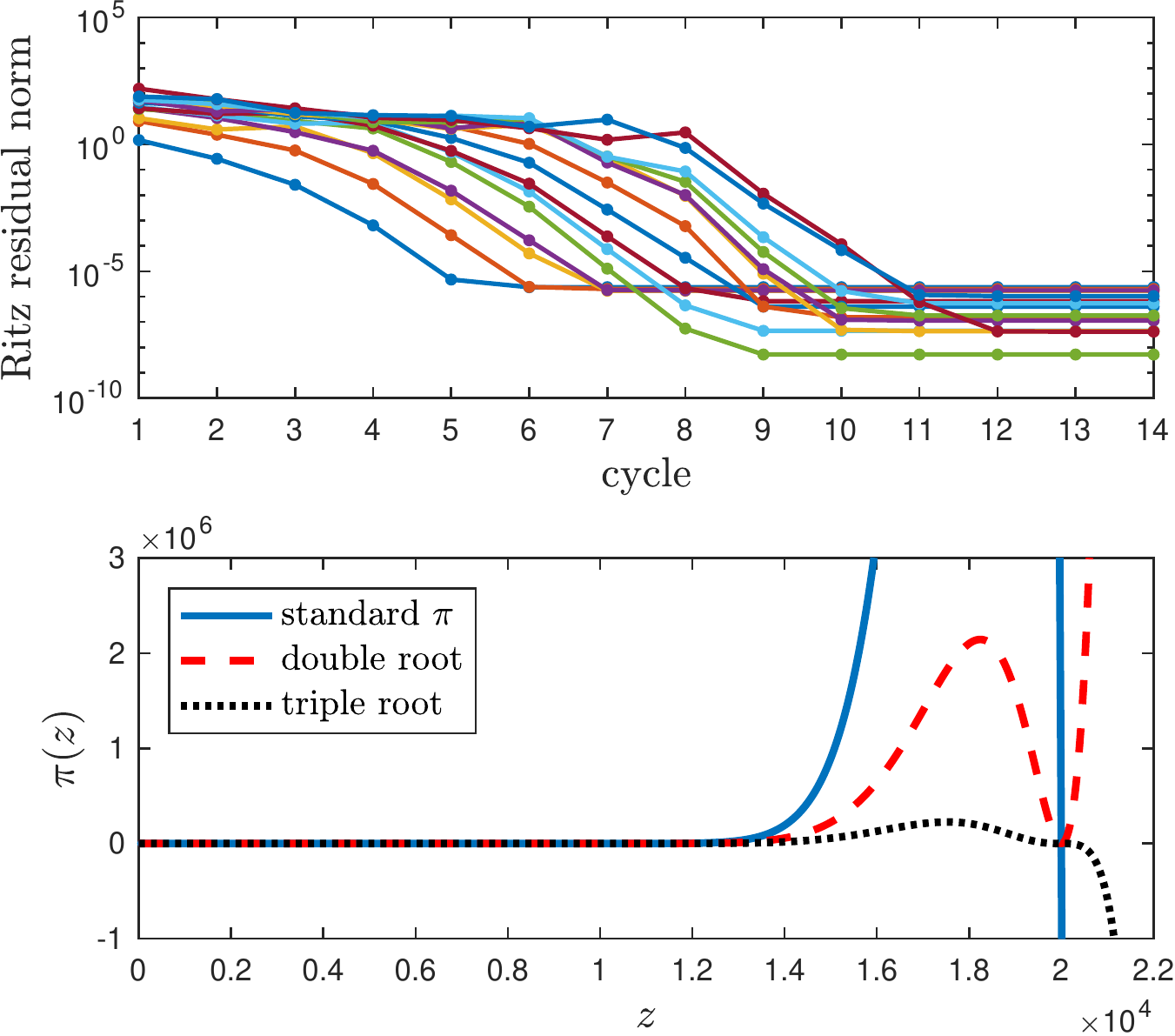}
\end{center}

\vspace*{-5pt}
\caption{\label{fig:ppaf11b}
Example~7, with $\deg=15$:  The top plot shows the residual norms, cycle by cycle: the large eigenvalue in $A$ limits the attainable accuracy.  The bottom plot shows how adding extra roots at the largest harmonic Ritz value tames the polynomial.}
%\vspace*{-5pt}
\end{figure}

%%%%%%%%%%%%%%%%%%%%%%%%%%%%%%%%%%%%%%%%%%%%%%%%%%%%%%%%%%%%%%%%%%%%%%%%%%%%%%%%%%%%%%%%%%%%%%%%%%
{\it Example 7.}  Let $A$ be the diagonal matrix of order $n=\mbox{10,000}$ with diagonal entries $0.1, 0.2, 0.3, \ldots, 9.9, 10, 11,$ $12, \ldots, 9908, 9909, 20000$, giving 100~relatively small eigenvalues and  one outlying eigenvalue, $\lambda=\mbox{20,000}$.
Using Arnoldi(50,20) with $\deg=5$, the $\nev=15$ computed eigenvalues all reach a residual norm at or below $\rtol=2.1\times 10^{-11}$, marginally better than the $6.1\times10^{-11}$ obtained without polynomial preconditioning.
With  $\deg=10$, the accuracy degrades to $7\times 10^{-9}$, while $\deg=15$ only reaches $2.3\times 10^{-6}$.
Figure~\ref{fig:ppaf11b} shows the $\deg=15$ residual convergence (top), and the corresponding preconditioner $\mrpoly$ (bottom, solid line): $\mrpoly$ has a root at $\theta_{15}=\mbox{20,000.0000000000036379}$, which of course is very near the large eigenvalue $\lambda_{10000}=\mbox{20,000}$.  When $\mrpoly(A)v$ is computed for some vector $v$, the factored form $\mrpoly(A) = \Pi_{j=1}^{15}(I - A/\theta_j)$ is used.  The component of $\mrpoly(A)v$  in the direction of the eigenvector $z_{10000}$ that corresponds to the large eigenvalue is $\gamma_{10000} \Pi_{i=1}^{15}(1 - \lambda_{10000}/\theta_i) z_{10000}$, where $\gamma_{10000}$ is the coefficient for $z_{10000}$ in an eigenvector expansion of $v$.  Fourteen of the $(1 - \lambda_{10000}/\theta_i)$ terms magnify this component and the fifteenth reduces it back down, but with substantial cancellation error: indeed, in this case $1-\lambda_{10000}/\theta_{15}\approx 2.22\times 10^{-16}$ (machine epsilon).

\begin{comment}
 degree d= 0; minimum of max residual norm = 6.13893e-11
 degree d= 5; minimum of max residual norm = 2.10731e-11
 degree d=10; minimum of max residual norm = 7.00051e-09
 degree d=15; minimum of max residual norm = 2.34915e-06
\end{comment}

%%%%%%%%%%%%%%%%%%%%%%%%%%%%%%%%%%%%%%%%%%%%%%%%%%%%%%%%%%%%%%%%%%%%%%%%%%%%%%%%%%%%%%%%%%%%%%%%%%
%IMPORTANT:  This figure is generated by PPArnoldiRootsStabExtraRoots program (copying roots to another program does not work!).  Wait, now have renamed the version of this that will create the figure as ppaf11B.
%\hspace{-.7in}
%\includegraphics[scale=.75]{ppaf11b}
%\vspace{-2.8in}
%%%%%%%%%%%%%%%%%%%%%%%%%%%%%%%%%%%%%%%%%%%%%%%%%%%%%%%%%%%%%%%%%%%%%%%%%%%%%%%%%%%%%%%%%%%%%%%%%%

We can monitor the possible loss of accuracy due to this cancellation error by computing how large the component is blown up by the other terms.  Define
\[
\pof(j) \equiv \prod_{\substack{i=1 \\ i \ne j }}^{d}\ \ |1 - \theta_{j}/\theta_i|;
\]
 ``$\pof@$'' stands for ``product of other factors evaluated at Ritz value.''
(This definition uses $\theta_j$ where we might like to use the $j$th desired eigenvalue, but $\theta_j$ will approximate it in the case of interest.)
The quantity $\pof(j)$ gauges the slope of $\pi$ at $\theta_j$, since
\[ |\pi'(\theta_j)| = \pof(j) / |\theta_j|.\]
(Unlike $\pi'(\theta_j)$, $\pof(j)$ is scale-invariant.)
Large $\pof(j)$ values signal points where $\mrpoly$ changes rapidly, warning of ill-conditioning in related computations.

\begin{figure}
%\vspace{-2.6in}
%\hspace{-.7in}
%\includegraphics[scale=.75]{ppaf11a}
%\vspace{-2.8in}
\begin{center}
\includegraphics[scale=.63]{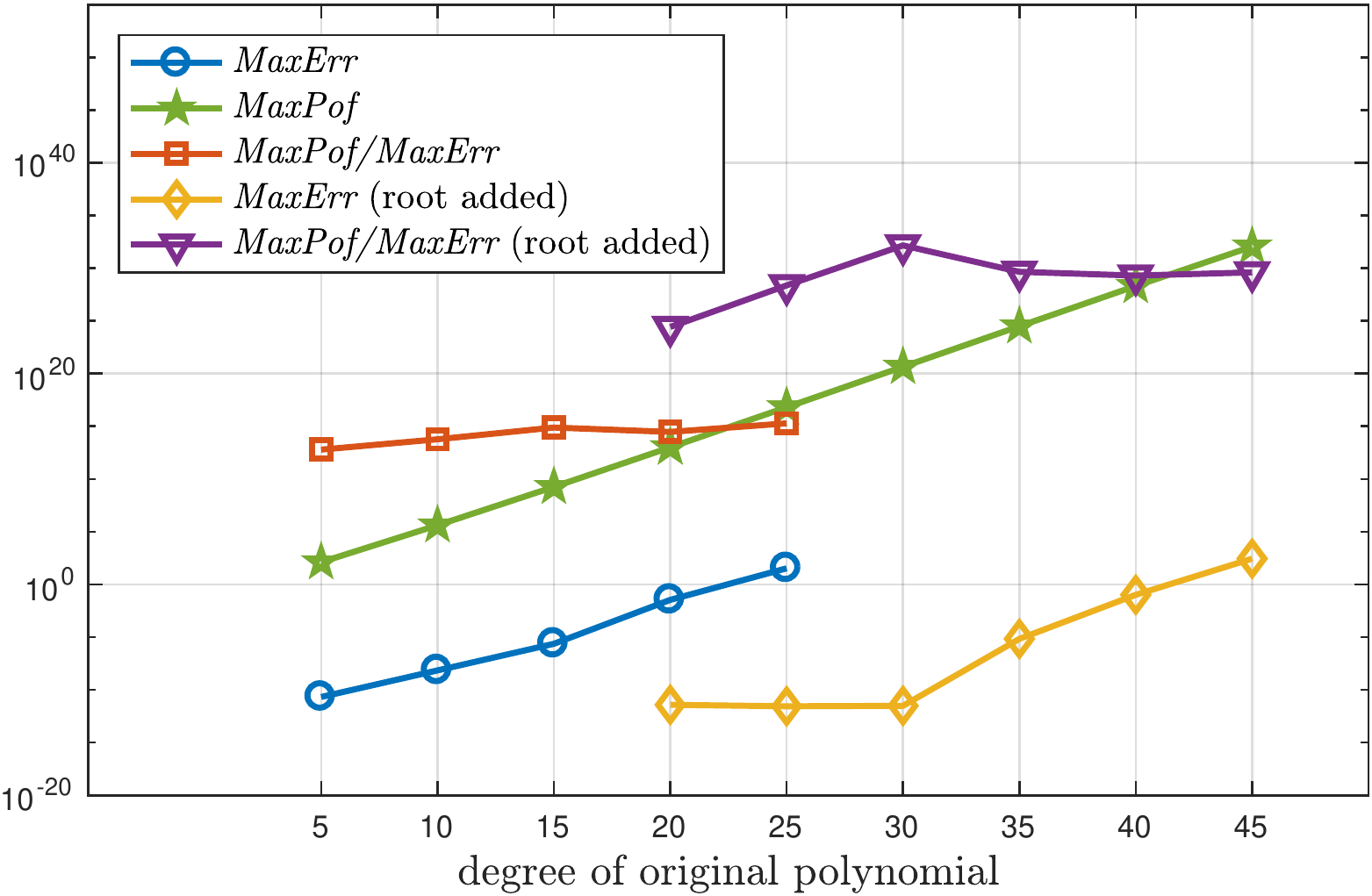}
\end{center}

\vspace*{-5pt}
\caption{\label{fig:ppaf11a}
Example~7, test of instability:  A large outlying eigenvalue gives a large $\pof(j)$; the attainable accuracy (blue circles) degrades as the degree $\deg$ increases.  An extra root to $\mrpoly$ at the largest harmonic Ritz value can improve the attainable accuracy for larger $\deg$ values (yellow diamonds).}
\vspace*{-15pt}
\end{figure}

Figure~\ref{fig:ppaf11a} shows that as the degree $\deg$ increases, the accuracy steadily degrades.  
For this example, the maximum Ritz residual norm grows with the maximum $\pof(j)$ value.  Let $\MaxErr$ be the maximum eventual residual norm of the 15~computed eigenvalues and let $\MaxPof$ be the maximum $\pof(j)$ value.  For the matrix in Example~7, Figure~\ref{fig:ppaf11a} plots $\MaxPof$ (which occurs at the largest harmonic Ritz value) with stars: it steadily increases as $\deg$ increases.  
Once the instability is the main source of error, starting at degree~10, the ratio $\MaxPof/\MaxErr$ is on the order of $10^{15}$. 

%%%%%%%%%%%%%%%%%%%%%%%%%%%%%%%%%%%%%%%%%%%%%%%%%%%%%%%%%%%%%%%%%%%%%%%%%%%%%%%%%%%%%%%%%%%%%%%%%%
%%%%%%%%%%%%%%%%%%%%%%%%%%%%%%%%%%%%%%%%%%%%%%%%%%%%%%%%%%%%%%%%%%%%%%%%%%%%%%%%%%%%%%%%%%%%%%%%%%

To improve stability when some $\pof(j)$ is large, we add an additional root at $\theta_j$ to $\mrpoly$, making $\theta_j$ a double root. When $\theta_j$ is almost an eigenvalue $\lambda$ of $A$, this makes $\mrpoly(\lambda)$ so near zero that even if the component of $\mrpoly(A)v$ in the direction of this eigenvector is off by several orders of magnitude, it is not significant relative to the other terms.  If $\theta_j$ is a double root of $\mrpoly$, then the slope $\mrpoly'(\theta_j)=0$, suggesting better conditioning.  However, the extra root increases the degree of $\mrpoly$ and the number of matrix-vector products with $A$ needed to apply $\mrpoly(A)$.  When is an extra root worth adding, and how many should be included?  In the test described above, an extra root added little benefit when $\pof(j)\le 10^4$.  In further testing, we have found that a new root is roughly needed for every factor of $10^{14}$ that $\MaxPof$ exceeds $10^4$.  Thus we suggest making $\theta_j$ a double root if $\pof(j)$ exceeds $10^4$, a triple root if $\pof(j)$ exceeds $10^{18}$, a quadruple root if $\pof(j)$ exceeds $10^{32}$, etc., incrementing by $10^{14}$ each time.

\vspace{.07in}
\begin{center}
\textbf{Adding Roots for Stability}
\end{center}
\begin{enumerate}
 \item {\bf Setup:} Assume the $d$ roots ($\theta_1, \ldots, \theta_d$) of $\mrpoly$ have been computed and then sorted according to the modified Leja ordering~\cite[alg.~3.1]{BaHuRe}.
 \item {\bf Compute $\pof(j)$:} For $j=1,\ldots,d$, compute $\pof(j) = \prod_{i \neq j} |1-\theta_j/\theta_i|$.  
 \item {\bf Add roots:} Compute least integer greater than $(\log_{10}(\pof(j)) - 4)/14$, for each $j$.  Add that number of $\theta_j$ values to the list of roots.  We add the first to the end of the list and if there are others, they are spaced into the interior of the current list, evenly between the occurrence of that root and the end of the list (keeping complex roots together).
\end{enumerate}
\vspace{.08in}

Other choices are possible for placing the new roots into the list of roots.  We also tried a second Leja sorting with the distance between identical roots defined to be a small amount such as $\alpha|\theta_j|$ for  $\alpha =10^{-15}$.  This was sometimes better and sometimes worse than the approach listed above.  This topic could use further investigation.

{\it Example 7 (continued) }   We now apply the procedure just given to Example~7, for increasing values of $\deg$.  The test adds a root all $\deg\ge 8$.  Figure~\ref{fig:ppaf11a} shows $\MaxErr$ for different degree polynomials with an added root (so the degree of the preconditioner is actually one more than the degree shown in the plot).  The accuracy for degree~25 is far better than without the added root ($2.8\times 10^{-12}$ compared to $3.4\times 10^{1}$).  However, even with a double root accuracy is lost for $\deg>30$; for large $\deg$, $\MaxPof/\MaxErr$ is roughly $10^{30}$.  For $\deg=40$ with one root added, $\MaxErr = 1.0\times 10^{-1}$.  However, at that point $\MaxPof = 2.0\times 10^{28}$, so according to the plan given above another extra root is needed.  With this triple root, $\MaxErr$ improves vastly to $4.5\times 10^{-12}$.
The bottom of Figure~\ref{fig:ppaf11b} compares the original $\deg=15$ polynomial to those with one and two added roots at the large harmonic Ritz value.  The slope of the original $\mrpoly$ is  large at $\lambda=\mbox{20,000}$.  Adding a root causes the polynomial to level off briefly there.  The polynomial with a triple root is not needed at this degree, but notice how it would add considerable stability near the extreme eigenvalue.

\begin{comment}

 MaxErr for d=25 (0 extra roots): 3.3759364e+01; maxpof = 6.4105973e+16
 MaxErr for d=25 (1 extra roots): 2.7709923e-12; maxpof = 6.4105973e+16
 MaxErr for d=40 (1 extra roots): 1.0351975e-01; maxpof = 2.0380107e+28
 MaxErr for d=40 (2 extra roots): 4.4706791e-12; maxpof = 2.0380107e+28
\end{comment}

We have tried this procedure of adding roots for realistic problems, including matrices used earlier in this paper and others, and it seems to work well.  The number of needed additional roots varies considerably.  For example, with the degree $\deg=100$, the Af23560 matrix from Example~2 uses 24~added roots (16~double and~4 triple), while the damped polynomial for S1rmq4m1 from Example~5 needs only 1~added root.  We do suspect there will be problems for which this procedure does not work; in fact, we have been able to devise such an example by skewing the starting vector against an outlying eigenvalue so severely that the associated harmonic Ritz value is far from that eigenvalue: extra roots in the wrong place do not much help.

\section{Double Polynomial Preconditioning}
%%%%%%%%%%%%%%%%%%%%%%%%%%%%%%%%%%%%%%%%%%%%%%%%%%%%%%%%%%%%%%%%%%%%%%%%%%%%%%%%%%%%%%%%%%%%%%%%%%

Communication-avoiding methods minimize operations that  transfer data across processors (and perform related synchronizations) such as dot products, potentially at the cost of  extra communication-free work on local processors.  We have already seen (e.g., Table~\ref{tbl:ppaf6}) that polynomial preconditioning can significantly reduce the number of dot products required to compute eigenvalues.  In fact, dot products can be even more significantly reduced by combining two levels of polynomial preconditioning, giving access to very high degree polynomials (which can permit lower subspace dimensions, and hence less memory). 

{\it Example 8.}  We revisit the convection-diffusion matrix from Example~1 of size $n = \mbox{640,000}$, again using Arnoldi(50,20) to compute the $\nev=15$ smallest eigenvalues.  
Table~\ref{tbl:cdo} reports results averaged over 10~trials, all of which converge toward all $\nev=15$ of the desired eigenvalues.)
Large degree preconditioning polynomials accelerate convergence, in terms of time and dot products.  
However, two concerns emerge as $\deg$ increases: construction of $\mrpoly$ becomes increasingly expensive (e.g., the $\deg = 150$ computation takes 60,127.8 dot products, 27,941.2 of which come from the GMRES run used to construct $\mrpoly$), and larger values of $\deg$ can cause stability problems.  We seek to avoid these limitations, while still reaping the benefits of high degree polynomials.

% d=150: dots = 60127.8; GMRES dots = 27941.2 

%%%%%%%%%%%%%%%%%%%%%%%%%%%%%%%%%%%%%%%%%%%%%%%%%%%%%%%%%%%%%%%%%%%%%%%%%%%%%%%%%%
\begin{table}[b!]

\caption{\label{tbl:cdo}
Example~8 (convection--diffusion, $n=\mbox{640,000}$).
Work required for Arnoldi(50,20) to compute the $\nev=15$ smallest magnitude eigenvalues.
Double polynomial preconditioning (bottom part of table) can significantly reduce the number of communication-intensive dot products required for standard polynomial preconditioning (top part of table).}

\begin{comment}
\begin{center}
\begin{tabular}{|c|c|c|c|c|c|c|}  \hline\hline

 degree  & cycles 	& $\mvps$  & $\cost$	& time	  &  dot products \\
 		 &  		& (thou's) & (thou's) & 	  &  (thousands)  \\ \hline
  \multicolumn{6}{|c|}{\emph{Polynomial Preconditioned Arnoldi}} \\ \hline
  1 (no prec) & 16,981 & 526	& 90,391	& 2.64 days		& 37,222 \\ \hline 
  10    &  471 	& 142  		& 3360 		& 2.01 hours	& 1033 \\ \hline
  25    &  113  & 85.8  	& 1120   	& 34.1 minutes	& 249  \\ \hline
  50    &  50   & 76.7  	& 730   	& 19.6 minutes	& 112  \\ \hline
  100   &  24   & 74.8  	& 588 		& 13.5 minutes	& 58.5 \\ \hline
  200   &  16   & 101  		& 735   	& 15.3 minutes	& 56.1 \\ \hline
  300   &  10   & 98.1  	& 733   	& 16.4 minutes	& 68.1 \\ \hline
  \multicolumn{6}{|c|}{\emph{Double Polynomial Preconditioning}}     \\ \hline
   $15\times20 = 300$  & 11    & 105 	& 692	& 13.4 minutes	& 24.9 \\ \hline
   $20\times25 = 500$  & 7     & 116   & 736   & 14.1 minutes  & 16.4 \\ \hline
   $20\times40 = 800$  & 4     & 113	& 702   & 13.2 minutes	& 10.3 \\ \hline
   $30\times55 = 1650$ & 2     & 134	& 819   & 14.9 minutes	& 6.91 \\ \hline
   $40\times65 = 2600$ & 1     & 139   & 809   & 14.0 minutes	& 5.69 \\ \hline
   \hline

\end{tabular}
\end{center}
\end{comment}

\vspace*{-3pt}
\begin{center}
\begin{tabular}{|r|c|c|c|c|c|c|}  \hline\hline
 \multicolumn{1}{|c|}{degree}  & cycles 	& $\mvps$  & $\cost$	& time	  &  dot products \\[-2pt]
 \multicolumn{1}{|c|}{$\deg$ or $d_1\times d_2$}		 &  		& (thousands) & (thousands) & (minutes)  &  (thousands)  \\ \hline
  \multicolumn{6}{|c|}{\emph{Polynomial Preconditioned Arnoldi}} \\ \hline
               $  0$ \phantom{0000} &   7436.1    &    223.1  &  42534.1  &    246.3  &  17178.2 \\ \hline 
               $ 10$ \phantom{0000}  &    259.0    &     78.3  &   1876.8  &     18.7  &    574.6 \\ \hline 
               $ 25$ \phantom{0000} &     84.3    &     64.2  &    845.4  &     11.6  &    188.6 \\ \hline 
               $ 50$ \phantom{0000} &     41.2    &     63.6  &    612.2  &     10.0  &     95.4 \\ \hline 
               $100$ \phantom{0000} &     20.6    &     64.8  &    524.1  &      9.5  &     57.9 \\ \hline 
               $125$ \phantom{0000} &     16.5    &     65.3  &    519.9  &      9.8  &     56.3 \\ \hline 
               $150$ \phantom{0000} &     14.0    &     67.0  &    535.0  &      9.9  &     60.1 \\ \hline 
                 \multicolumn{6}{|c|}{\emph{Double Polynomial Preconditioning}}     \\ \hline
 $ 15\times 20= \phantom{0}300$  &      3.8    &     41.0  &    273.6  &      1.7  &      9.7 \\ \hline 
 % $ 15\times 30= \phantom{0}450$   &      2.5    &     43.4  &    282.4  &      1.8  &      7.3 \\ \hline 
 $ 15\times 40= \phantom{0}600$   &      2.0    &     48.9  &    314.8  &      2.0  &      6.9 \\ \hline 
 % $ 25\times 30= \phantom{0}750$   &      2.0    &     60.9  &    384.6  &      2.5  &      6.6 \\ \hline 
 $ 15\times 50= \phantom{0}750$   &      2.0    &     61.2  &    392.0  &      2.5  &      7.9 \\ \hline 
 % $ 15\times 60= \phantom{0}900$   &      1.0    &     45.9  &    295.6  &      1.9  &      6.8 \\ \hline 
 $ 25\times 40=1000$   &      1.0    &     51.0  &    321.0  &      2.1  &      5.2 \\ \hline 
 % $ 25\times 50=1250$   &      1.0    &     63.8  &    400.9  &      2.5  &      6.5 \\ \hline 
 $ 25\times 60=1500$   &      1.0    &     76.5  &    481.4  &      3.1  &      8.0 \\ \hline 
 % $ 35\times 50=1750$   &      1.0    &     89.3  &    556.0  &      3.5  &      7.5 \\ \hline 
 \hline
\end{tabular}
\end{center}
\vspace*{-15pt}
\end{table}
%%%%%%%%%%%%%%%%%%%%%%%%%%%%%%%%%%%%%%%%%%%%%%%%%%%%%%%%%%%%%%%%%%%%%%%%%%%%%%%%%%

These observations motivate \emph{Double Polynomial Preconditioning}.  Start by generating the GMRES polynomial $\mrpoly_1$ of degree $\deg_1$ for $A$.\ \ As before, we expect the smallest magnitude eigenvalues of $A$ to be mapped to the eigenvalues of $\mrpoly_1(A)$ nearest~1.  Thus define $\tau(\alpha) \equiv 1 - \mrpoly_1(\alpha)$:  we seek the smallest magnitude eigenvalues of $\tau(A)$.
To compute these smallest magnitude eigenvalues of $\tau(A)$, apply polynomial preconditioning to this matrix, i.e., apply GMRES to $\tau(A)$ to generate a new GMRES polynomial $\mrpoly_2$ of degree $\deg_2$.  Now apply Arnoldi$(m,k)$ to compute the $\nev$ eigenvalues of $\mrpoly_2(\tau(A))$ nearest~1.
The composite polynomial $\mrpoly_2(\tau(\cdot))$ has degree $d_1@d_2$, making use of extremely high degree polynomials more practical.

{\it Example 8 (continued).} The bottom half of Table~\ref{tbl:cdo} shows the effectiveness of double polynomial preconditioning for the convection-diffusion problem.  The first column reports the polynomial degrees $d_1$ and $d_2$; e.g., $15\times 20 = 300$ means that $\tau$ has degree $\deg_1 = 15$ and $\pi_2$ has degree $\deg_2=20$, so the composite polynomial $\pi_2(\tau(\cdot))$ has degree~300. Because high degree composite polynomials can be formed without the need for much GMRES orthogonalization, the dot products are greatly reduced (but other costs can go up).  
For composite degree~$25\times40=1000$, only  5,231.2~dot products are needed, a ten-fold reduction from the lowest number given for single polynomial preconditioning.  Arnoldi(50,20) needs only one cycle with this high degree polynomial.  (In these tests, we only check residual norms at the end of cycles.  To further reduce operations, we could check residuals mid-cycle and terminate early.)

%{\it Example 11.5.} ((If we decide to add what is now Table 8.2 ))

{\it Example 9.} 
We revisit the matrix in Example~7.  In this case, double polynomial preconditioning can help cure the instability, though this is not guaranteed in general. Let the first polynomial have degree $d_1 = 6$, which has a root near $\mbox{19,991.2}$, close enough to $\lambda=\mbox{20,000}$ so that the spectrum of $\pi_2(\tau(A))$ with $d_2=20$ does not have an outstanding eigenvalue, and there is no instability: Arnoldi(50,20) finds the $\nev=15$ smallest eigenvalues in two~cycles (composite degree $6\times 20 = 120$). Next, we change to $d_1 = 5$, for which $\pi_1$ has a root only at~19,700.5: not close enough to $\lambda=\mbox{20,000}$, so for $d_2=20$ the matrix $\pi_2(\tau(A))$ has an outstanding eigenvalue, and no progress is made in 50~Arnoldi cycles.  However, the $\MaxPof$ test described in Section~7 suggests that a double root be added here:  that is sufficient to give convergence in two cycles, finding most of the $\nev=15$ desired eigenvalues.

Further investigation is needed of stability for double polynomial preconditioning and whether the same test we applied for one polynomial remains effective here.

\section{Conclusions}
%%%%%%%%%%%%%%%%%%%%%%%%%%%%%%%%%%%%%%%%%%%%%%%%%%%%%%%%%%%%%%%%%%%%%%%%%%%%%%%%%%%%%%%%%%%%%%%%%%

Polynomial preconditioning can vastly improve eigenvalue calculations for difficult problems, giving the benefit of working with high-degree polynomials in $A$ without requiring high-dimensional subspaces.  

We have focused on preconditioning with the GMRES (residual) polynomial, which is easy to compute and adapts to the spectrum of $A$.  Our computational experiments illustrated the success of this method and identified a few scenarios that can be remedied with special handling:  using multiple starting vectors for GMRES; damping the GMRES starting vector; adding extra copies of certain roots to enhance stability.
While polynomial preconditioning often reduces matrix-vector products for difficult problems, the reduction in vector operations such as dot products is even greater, so this approach holds great promise for communication-avoiding eigenvalue computation on high performance computers.  \emph{Double Polynomial Preconditioning} gives access to high degree polynomials in $A$, and can further reduce dot products.  Techniques from~\cite{AvoidComm,Ho10} can potentially aid parallel implementations. 

Future research should include computation of interior eigenvalues, generalized eigenvalue problems, and application to computing stable eigenvalues in matrices that exhibit a significant departure from normality~\cite[chap.~28]{TE05}.   Stability control for the double polynomial preconditioning should also be investigated.   

\begin{comment}
\section*{Acknowledgments} The authors appreciate . . .

\end{comment}


\begin{thebibliography}{10}

\bibitem{BaHuRe}
{\sc Z.~Bai, D.~Hu, and L.~Reichel}, {\em A {N}ewton basis {GMRES}
  implementation}, IMA J. Numer. Anal., 14 (1994), pp.~563--581.

\bibitem{BER04}
{\sc C.~Beattie, M.~Embree, and J.~Rossi}, {\em Convergence of restarted
  {Krylov} subspaces to invariant subspaces}, SIAM J.~Matrix Anal.\ Appl., 25
  (2004), pp.~1074--1109.

\bibitem{BES05}
{\sc C.~A. Beattie, M.~Embree, and D.~C. Sorensen}, {\em Convergence of
  polynomial restart {Krylov} methods for eigenvalue computations}, SIAM
  Review, 47 (2005), pp.~492--515.

\bibitem{Cha83}
{\sc F.~Chatelin}, {\em Spectral Approximation of Linear Operators}, Academic
  Press, New York, 1983.

\bibitem{AvoidComm}
{\sc J.~Demmel, M.~Hoemmen, M.~Mohiyuddin, and K.~Yelick}, {\em Avoiding
  communication in sparse matrix computations}, in 2008 IEEE International
  Symposium on Parallel and Distributed Processing, IEEE, 2008.

\bibitem{DM12}
{\sc J.~Duintjer~Tebbens and G.~Meurant}, {\em Any {Ritz} value behavior is
  possible for {Arnoldi} and {GMRES}}, SIAM J.~Matrix Anal.\ Appl., 33 (2012),
  pp.~958--978.

\bibitem{Emb09}
{\sc M.~Embree}, {\em The {Arnoldi} eigenvalue iteration with exact shifts can
  fail}, SIAM J.~Matrix Anal.\ Appl., 31 (2009), pp.~1--10.

\bibitem{Fr92}
{\sc R.~W. Freund}, {\em Quasi-kernel polynomials and their use in
  non-{H}ermitian matrix iterations}, J. Comput. Appl. Math., 43 (1992),
  pp.~135--158.

\bibitem{Ho10}
{\sc M.~F. Hoemmen}, {\em Communication-avoiding {K}rylov subspace methods}.
\newblock PhD Thesis, EECS Dept., University of California at Berkeley, 2010.

\bibitem{Kat80}
{\sc T.~Kato}, {\em Perturbation Theory for Linear Operators}, Springer-Verlag,
  Berlin, corrected second~ed., 1980.

\bibitem{La52B}
{\sc C.~Lanczos}, {\em Chebyshev polynomials in the solution large-scale linear
  systems}, Proc. ACM,  (1952), pp.~124--133.

\bibitem{LiXiVeYaSa}
{\sc R.~Li, Y.~Xi, E.~Vecharynski, C.~Yang, and Y.~Saad}, {\em A thick-restart
  {L}anczos algorithm with polynomial filtering for {H}ermitian eigenvalue
  problems}, SIAM J. Sci. Comput., 38 (2016), pp.~A2512--A2534.

\bibitem{PPG}
{\sc Q.~Liu, R.~B. Morgan, and W.~Wilcox}, {\em Polynomial preconditioned
  {GMRES} and {GMRES-DR}}, SIAM J. Sci. Comput., 37 (2015), pp.~S407--S428.

\bibitem{MSR94}
{\sc K.~Meerbergen, A.~Spence, and D.~Roose}, {\em Shift-invert and {Cayley}
  transforms for detection of rightmost eigenvalues of nonsymmetric matrices},
  BIT, 34 (1994), pp.~409--423.

\bibitem{IE}
{\sc R.~B. Morgan}, {\em Computing interior eigenvalues of large matrices},
  Linear Algebra Appl., 154-156 (1991), pp.~289--309.

\bibitem{Arnoldi-R}
\leavevmode\vrule height 2pt depth -1.6pt width 23pt, {\em On restarting the
  {A}rnoldi method for large nonsymmetric eigenvalue problems}, Math. Comp., 65
  (1996), pp.~1213--1230.

\bibitem{HRAM}
{\sc R.~B. Morgan and M.~Zeng}, {\em A harmonic restarted {A}rnoldi algorithm
  for calculating eigenvalues and determining multiplicity}, Linear Algebra
  Appl., 415 (2006), pp.~96--113.

\bibitem{NaReTr}
{\sc N.~M. Nachtigal, L.~Reichel, and L.~N. Trefethen}, {\em A hybrid {GMRES}
  algorithm for nonsymmetric linear systems}, SIAM J. Matrix Anal. Appl., 13
  (1992), pp.~796--825.

\bibitem{PaPavdV}
{\sc C.~C. Paige, B.~N. Parlett, and H.~A. van~der Vorst}, {\em Approximate
  solutions and eigenvalue bounds from {K}rylov subspaces}, Numer. Linear
  Algebra Appl., 2 (1995), pp.~115--133.

\bibitem{PaSa}
{\sc C.~C. Paige and M.~A. Saunders}, {\em Solution of sparse indefinite
  systems of linear equations}, SIAM J. Numer. Anal., 12 (1975), pp.~617--629.

\bibitem{Rutis}
{\sc H.~Rutishauser}, {\em Theory of gradient methods}, in Refined Iterative
  Methods for Computation of the Solution and the Eigenvalues of Self-Adjoint
  Boundary Value Problems, M.~Engeli, T.~Ginsburg, H.~Rutishauser, and
  E.~Stiefel, eds., Birkhauser, Basel, 1959, pp.~24--49.

\bibitem{Saa80}
{\sc Y.~Saad}, {\em Variations on {Arnoldi's} method for computing
  eigenelements of large unsymmetric matrices}, Linear Algebra Appl., 34
  (1980), pp.~269--295.

\bibitem{Sa84b}
{\sc Y.~Saad}, {\em {C}hebychev acceleration techniques for solving large
  nonsymmetric eigenvalue problems}, Math. Comp., 42 (1984), pp.~567--588.

\bibitem{Sa87b}
\leavevmode\vrule height 2pt depth -1.6pt width 23pt, {\em Least squares
  polynomials in the complex plane and their use for solving sparse
  nonsymmetric linear systems}, SIAM J. Numer. Anal., 24 (1987), pp.~155--169.

\bibitem{Saa03}
{\sc Y.~Saad}, {\em Iterative Methods for Sparse Linear Systems}, SIAM,
  Philadelphia, second~ed., 2003.

\bibitem{Sa11}
{\sc Y.~Saad}, {\em Numerical Methods for Large Eigenvalue Problems, 2nd
  Edition}, SIAM, Philadelphia, PA, 2011.

\bibitem{SaSc}
{\sc Y.~Saad and M.~H. Schultz}, {\em {GMRES}: a generalized minimum residual
  algorithm for solving nonsymmetric linear systems}, SIAM J. Sci. Stat.
  Comput., 7 (1986), pp.~856--869.

\bibitem{So}
{\sc D.~C. Sorensen}, {\em Implicit application of polynomial filters in a
  $k$-step {A}rnoldi method}, SIAM J. Matrix Anal. Appl., 13 (1992),
  pp.~357--385.

\bibitem{StSaWu}
{\sc A.~Stathopoulos, Y.~Saad, and K.~Wu}, {\em Dynamic thick restarting of the
  {D}avidson, and the implicitly restarted {A}rnoldi methods}, SIAM J. Sci.
  Comput., 19 (1998), pp.~227--245.

\bibitem{St01}
{\sc G.~W. Stewart}, {\em A {K}rylov--{S}chur algorithm for large
  eigenproblems}, SIAM J. Matrix Anal. Appl., 23 (2001), pp.~601 -- 614.

\bibitem{Sti58}
{\sc E.~L. Stieffel}, {\em Kernel polynomials in linear algebra and their
  numerical applications}, U. S. Nat. Bur. Standards, Appl. Math. Ser., 49
  (1958), pp.~1--22.

\bibitem{Tho06}
{\sc H.~K. Thornquist}, {\em Fixed-Polynomial Approximate Spectral
  Transformations for Preconditioning the Eigenvalue Problem}, PhD thesis, Rice
  University, 2006.
\newblock Technical report TR06-05, Department of Computational and Applied
  Mathematics, Rice University.

\bibitem{TE05}
{\sc L.~N. Trefethen and M.~Embree}, {\em Spectra and Pseudospectra: The
  Behavior of Nonnormal Matrices and Operators}, Princeton University Press,
  Princeton, NJ, 2005.

\bibitem{WuSi}
{\sc K.~Wu and H.~Simon}, {\em Thick-restart {L}anczos method for symmetric
  eigenvalue problems}, SIAM J. Matrix Anal. Appl., 22 (2000), pp.~602 -- 616.

\end{thebibliography}
\end{document}